\documentclass{article}
\usepackage[]{epsf,epsfig,amsmath,amssymb,amsfonts,latexsym}
\title{Volumetric Properties of the Convex Hull of an $n$-dimensional Brownian Motion}
\author{Ronen Eldan \thanks{Partially supported by the Israel Science Foundation}}
\newtheorem{theorem}{Theorem}
\newtheorem{lemma}[theorem]{Lemma}
\newtheorem{question}[theorem]{Question}

\newtheorem{corollary}[theorem]{Corollary}

\newtheorem{remark}[theorem]{Remark}

\def \PP {\mathbb P}
\def\qed{\hfill $\vcenter{\hrule height .3mm
\hbox {\vrule width .3mm height 2.1mm \kern 2mm \vrule width .3mm
height 2.1mm} \hrule height .3mm}$ \bigskip}
\def\P{\mathbb{P}}

\def\EE{\mathbb{E}}

\def\R{\mathbb{R}}
\def\RR{\mathbb{R}}

\def\Sph{S^{n-1}}

\def\Vol{\mathrm{Vol}}
\def\Conv{\mathrm{Conv}}
\def\BBN{\mathbb{B}_n}

\begin{document}
\maketitle

\begin{abstract}
Let $K$ be the convex hull of the path of a standard Brownian motion $B(t)$ in $\RR^n$, taken at time $0 \leq t \leq 1$. We derive formulas for the expected volume and surface area of $K$. Moreover, we show that in order to approximate $K$ by a discrete version of $K$, namely by the convex hull of a random walk attained by taking $B(t_n)$ at discrete (random) times, the number of steps that one should take in order for the volume of the difference to be relatively small is of order $n^3$. Next, we show that the distribution of facets of $K$ is in some sense scale invariant: for any given family of simplices (satisfying some compactness condition), one expects to find in this family a constant number of facets of $t K$ as $t \to \infty$. Finally, we discuss some possible extensions of our methods and suggest some further research.
\end{abstract}
\bigskip
\section{Introduction}

Convex hulls seem to attract a significant amount of interest, in some cases for representing physical phenomena and in others for their central importance in many algorithmic methods. Random convex hulls, in numerous different settings, have been widely studied by probabilists and geometers (see \cite{Bar,CMR,Reizner,S} for surveys of the subject). One example of a random convex hull that has been studied is the convex hull of a Brownian motion, which may represent, for example, the domain of influence on a diffusing particle in a certain physical system. The object of this paper is to further study the object generated by taking the convex hull of the standard Brownian motion in $\RR^n$. \\

The convex hull of the path of the planar Brownian motion has been quite extensively studied. Much is known about this object, including its expected area and perimeter length, the degree of smoothness of its boundary, the rate of convergence of the area of a convex hull of a random walk to its area, etc. (see e.g. \cite{CMR,Bax,BL,T,Elb,CHM} and references therein). However, it seems like much less is known about the convex hull of the Brownian motion in higher dimensions. Two examples of notable works concerning the higher dimensional case are a paper by Kampf, Last and Molchanov, \cite{KLM} in which, for instance, the first and second intrinsic volumes are calculated and a work by Kinney, \cite{K}, in which a bound for the total curvature is established. \\

We extend some known results from the planar case to the higher dimensional case, as well as obtain certain asymptotics of the behaviour of these objects as the dimension goes to infinity. We introduce new methods which may be further used to study volumetric and combinatorial properties of the convex hull of the Brownian motion and random walk. \\

Since not many concrete examples of high-dimensional convex bodies are known, finding new explicit constructions may lead to a deeper understanding of the theory around such bodies and may possibly provide counter examples to some general conjectures related to those bodies (one example of such a conjecture is the \emph{hyperplane conjecture} described in \cite{kk}). The methods introduced in the present paper may be seen as a basic tool box for the geometer to analyze the convex hull of a high dimensional Brownian motion, thus enhancing the understanding of this specific explicit construction. As explained in the discussion below, one may expect this construction to have some properties that are different from most known examples, thus making it an interesting case to study. For instance, there is evidence which suggests that it admits highly-neighborly approximating polytopes, however its boundary is smooth. See Section \ref{secdiscussion} for more details and examples. \\

Let us introduce our setting. Fix a dimension $n \in \mathbb{N}$. Let $B(\cdot)$ be a standard Brownian motion in $\RR^n$. For a subset $A \subset \RR^n$, by $\Conv(A)$ we denote the minimal convex set containing $A$, the \emph{convex hull} of $A$. Our main object of concern in this paper will be
$$
K = \Conv(\{B(t), 0 \leq t \leq 1 \}).
$$
We will be interested, for example, in its expected volume and in its distribution of facets. In order to better study these properties, in many cases it will be convenient to introduce an approximation for this object by a simpler object, namely the convex hull of a random walk. We construct the random walk as follows: \\ \\
Let $P = ((x_1,y_1), (x_2,y_2),...)$ be a Poisson point process of intensity $1$ in the set $[0,1] \times [0,\infty]$ and for all $\alpha \geq 0$, define
$$
\Lambda_\alpha = \left \{ x_i| ~ y_i \leq \alpha, ~ i \in \mathbb{N} \right \} \cup \{0,1\}.
$$
The process $\Lambda$ can be thought of as a "Poisson rain" on the interval $[0,1]$: note that for all $\alpha \geq 0$, $\Lambda_\alpha$ 
is a Poisson point-process of intensity $\alpha$ on the unit interval and that the family $\Lambda_\alpha$ is increasing with $\alpha$. For a fixed value of $\alpha$, writing $\Lambda_\alpha = (t_1,...,t_N)$ where $0=t_1 \leq ... \leq t_N=1$, we can think of $(B(t_1), B(t_2), ..., B(t_N))$ as a random walk in $\RR^n$. Finally, for all $\alpha > 0$, we define
$$
K_\alpha = \Conv(\{B(t)| ~t \in \Lambda_\alpha \} ),
$$
so $K_\alpha$ is a monotone sequence of discrete approximations of $K$, each defined as the convex hull of a certain random walk. \\ \\

For a measurable set $L \subset \RR^n$, we denote the $k$-dimensional Hausdorff measure of $L$ by $\Vol_k(L)$. By $\partial L$ we denote the boundary of $L$. The first theorem we prove is a formula for the expected volume and surface area of $K$. The theorem reads,

\begin{theorem} \label{thmvolume}
One has, for every dimension $n \geq 2$,
$$
\EE[\Vol_n(K)] = \left (\frac{\pi}{2} \right )^{n/2} \frac{1}{\Gamma \left (\frac{n}{2} + 1 \right )^2}
$$
and
$$
\EE[\Vol_{n-1}(\partial K)] = \frac{2 (2 \pi)^{(n-1) / 2}}{\Gamma(n)}.
$$
\end{theorem}

As was pointed to us by Christoph Th\"{a}le, a corollary of this theorem is an explicit formula for expressing all intrinsic volumes of the body $K$, which generalizes a result of Kampf, Last and Molchanov (\cite{KLM}). Let $V_j(K)$ denote the $j$-th intrinsic volume of $K$ (see \cite{KLM} for a definition).

\begin{corollary} \label{cormixed}
One has for all $n \geq 2$ and for all $1 \leq j \leq n$
$$
\EE[V_j(K)] = \left ( n \atop j \right ) \left ( \frac{\pi}{2} \right )^{j/2} \frac{\Gamma \left (\frac{n-j}{2} + 1 \right ) }{ \Gamma \left (\frac{j}{2} + 1 \right ) \Gamma \left (\frac{n}{2} + 1 \right )}.
$$ 
\end{corollary}

Our next result is a derivation of asymptotics for the number of steps needed in order to approximate the convex hull of the Brownian motion, $K$, by the convex hull of the random walk, $K_\alpha$. Our theorem roughly states that the correct order of points needed in order for the volume of $K_\alpha$ to be a proportion of the volume of $K$ is $n^3$. It reads,

\begin{theorem} \label{thmapprox}
One has the following bounds: For all $n \geq 2$ and all $\alpha > 0$,
\begin{equation} \label{aproxlower}
\frac{\EE[\Vol_n(K \setminus K_\alpha)]}{\EE[\Vol_n(K)] } \leq e^{-n} + 16 \sqrt{\frac{n^{3}}{\alpha}}.
\end{equation}
On the other hand, for all $\alpha < n^3 / 8$, one has
\begin{equation} \label{aproxupper}
\frac{\EE[\Vol_n(K_\alpha)]}{\EE[\Vol_n(K)]} \leq 100 \frac{\alpha}{n^3} \log^2 \left (\frac{n^3}{\alpha} \right) . 
\end{equation}
\end{theorem}
Note that, according to the above theorem, for any given proportion constant $R<1$, there exists a constant $C(R)$ independent of the dimension,
such that whenever $\alpha > C(R) n^3$, the proportion between the expected volume $K_\alpha$ out of the entire volume of $K$ will be at least $R$. By basic properties of the Poisson process, the same will be true if we take $C(R) n^3$ uniform points on the interval $[0,1]$. On the other hand, the second part of the theorem shows us that taking only $o(n^3)$ points will yield $\EE[ \Vol_n(K_\alpha)] = o( \EE[ \Vol_n(K)])$. \\ \\

Our last result concerns with the distribution of facets of $K$. In order to formulate it, we need some notation. For two $(n-1)$-dimensional simplices $s_1, s_2 \subset \RR^n$ we say that $s_1$ and $s_2$ are equivalent if they are equal up to some translation. From this point further, by slight abuse of terminology, the term simplex will refer to an equivalence class of simplices. We denote by $\mathcal{S}$ the set of $(n-1)$-dimensional simplices and let $\mathcal{F}(K) \subset \mathcal{S}$ denote the set of $(n-1)$-dimensional facets of $K$ (i.e., the set of $(n-1)$-dimensional simplices lying entirely in the boundary of $K$, which are \emph{maximal} in the sense that they are not strictly contained in any simplex lying in the boundary of $K$). For a family of simplices $\mathcal{C} \subset \mathcal{S}$, we define
$$
M_K(\mathcal{C}) = \EE \left [ \# (\mathcal{F}(K) \cap \mathcal{C})  \right],
$$
the expected number of facets of $K$ which are in $\mathcal{C}$. Our aim is to study the behaviour of $M_K(\mathcal{C})$. \\ \\
Next, for a set $L \subset \RR^n$ and for $\epsilon > 0$, we denote
$$
e(L,\epsilon) := \{x \in \RR^n | ~\exists y \in L, ~|y - x| \leq \epsilon \},
$$
the $\epsilon$-extension of $L$. For two sets $L,T \subset \RR^n$, we denote
$$
d_H(L,T) := \inf \{\epsilon; ~ L \subset e(T,\epsilon) \mbox{ and } T \subset e(L,\epsilon) \},
$$
the Hausdorff distance between $L$ and $T$. For a family $\mathcal{C} \subset \mathcal{S}$ we say that $\mathcal{C}$ is \textbf{compact} if the set is compact with respect to the topology induced by the Hausdorff metric. We say that $\mathcal{C}$ is \textbf{non-degenerate} if there exists a constant $c>0$ such that every
simplex $s \in \mathcal{C}$ satisfies $\Vol_{n-1}(s) \geq c$. For a set $\mathcal{C} \subset \mathcal{S}$ and $t > 0$, we understand $t \mathcal{C}$ as $\{t s; ~ s \in \mathcal{C} \}$. \\ \\
The body $K$ is not a polytope; the set $\mathcal{F}(K)$ is almost-surely infinite. This follows directly, for example, from the result of \cite{CHM} which implies that
two-dimensional projections of $K$ almost-surely have a smooth boundary. In fact, we conjecture that $\partial K$ is smooth in any dimension (see discussion in Section \ref{secdiscussion}). However, we do know that up to a set of measure $0$, all of the $n-1$-dimensional volume of the boundary of $K$ is contained in the interior of these facets (see Corollary \ref{corfacets} below). Our last theorem, which characterizes the distribution of the boundary facets, reads

\begin{theorem} \label{thmscaling}
Let $\mathcal{C} \subset \mathcal{S}$ be a family of simplices. The function
$$
t \to M_K(t \mathcal{C})
$$
is decreasing. Moreover, if $\mathcal{C}$ is compact and non-degenerate, then the above function is bounded from above, and thus the limit
$$
\lim_{t \to 0^+} M_K(t \mathcal{C})
$$
exists.
\end{theorem}

Roughly speaking, the above theorem states that one should expect to find a constant number of facets of a given shape at any scale. For example, if
$n=3$ and the set $\mathcal{C}$ consists of all triangles of distance $\epsilon$ to an equilateral triangle whose edge has length $1$, the theorem suggests that there exists some constant $S$ such that one expects to find approximately one almost-equilateral facet whose edges are of length between $t$ and $tS$, for all small enough values of $t$. 
\\ \\

Some of our methods of proof extend a certain formula that appears in \cite{E1}, based on very simple principles from integral geometry. Some of these principles have already been used by Baxter \cite{Bax} in order to study the convex hull of planar random walks. The structure of this paper is the following: in Section 2, we derive certain estimates for one-dimensional random walks, which will be used later on. In Section 3 we establish some formulas concerning the distribution of facets of the polytope $K_\alpha$, which will be one of the central ingredients in our proofs. In Section 3 we prove Theorem \ref{thmvolume} and Corollary \ref{cormixed}. In sections 5 and 6 we prove theorems \ref{thmapprox} and \ref{thmscaling} respectively. Finally, in Section 7 we discuss some further possible extensions of our methods and raise some questions for further research. \\ \\

Throughout this paper, the symbols $C,C',C'',c,c',c''$ denote positive universal constants whose values may change between different formulas. Given a subset $A \subset \RR^n$, by $\Conv(A)$ we denote the convex hull of $A$, $\partial A$ will denote its boundary, $Cl(A)$ its closure and $Int(A)$, its interior. For a function $f : \RR^n \to \R$ we write $supp(f) = Cl(\{x; f(x) \neq 0 \})$, its support. By slight abuse of terminology, the word polytope in this paper refers to a convex polytope, i.e., the convex hull of the finite number of points.  \\

\section{One dimensional random walks}

In this section we derive some estimates concerning one-dimensional random walks. 

Let $0 \leq t_1 \leq ... \leq t_N \leq 1$ be a Poisson point process on $[0,1]$ with intensity $\alpha$ and let $B(t)$ be a standard 1-dimensional Brownian motion, independent from the above Poisson process. Consider the random walk  $B(0), B(t_1),...,B(t_N)$. By slight abuse of notation, for $1 \leq j \leq n$, denote $B(j) = B(t_j)$. 
Let us calculate the probability that $B(j) \geq 0$ for all $1 \leq j \leq N$. \\
Define a random variable,
$$
X = \int_0^1 \mathbf{1}_{\{ B(t) < 0 \}} dt.
$$
Recall the second arcsine law of P. L\'{e}vy (see for example \cite{MP}, Chapter 5, p. 137), according to which, $X$ has the distribution whose density $f_X$ satisfies 
$$
f_X(x) = \frac{1}{\pi} (x-x^2)^{-1/2} \mathbf{1}_{x \in [0,1]}.
$$
By definition of the Poisson distribution, we know that for all measurable $A \subset (0,1)$ one has
$$
\PP(\{t_1,...t_N\} \cap A = \emptyset ) = e^{-\alpha |A|}
$$
where $|A|$ denotes the Lebesgue measure of $A$. Applying this to $A = \{t \in (0,1);~  B(t) < 0 \} $ gives
$$
\P(B(t_i) \geq 0, ~~ \forall 1 \leq i \leq N) = \EE \left [e^{-\alpha X} \right ] = \frac{1}{\pi} \int_0^1  e^{-\alpha x} (x-x^2)^{-1/2} dx =
$$
(substituting $t^2 = \alpha x$)
$$
\frac{2}{\pi \sqrt{\alpha}} \int_0^{\sqrt{\alpha}} e^{-t^2} \frac{1}{\sqrt{1 - \frac{t^2}{\alpha}}} dt.
$$
$$
~
$$
Now suppose that $W(t)$ is a Brownian bridge such that $W(0)=W(1)=0$ and consider the discrete Brownian bridge 
$W(0), W(t_1),...,W(t_N), W(1)$. 

The cyclic shifting principle (see e.g., \cite{Bax}) is the following observation: for every $0 \leq s \leq 1$, define
$\Gamma_s(t) = t + s$, where the sum is to be understood as a sum on the torus $[0,1]$. Then the function $W \circ \Gamma_s (t) - W(s)$ has 
the same distribution as the function $W(t)$. Now, since there is exactly one choice $i$ between $0$ and $N$ such that
$W(t_j) - W(t_i)$ will be non-negative for every $1 \leq j \leq N$ (where $t_0 = 0$), it follows that for only one choice of $0 \leq i \leq N$, the function
$$
W \circ \Gamma_{t_i}(\cdot) - W(t_i)
$$
will be positive for all the points $t_j - t_i$, $0 \leq j \leq N$ (where the subtraction is again understood on the torus $[0,1]$). Since
the points $t_1,...,t_N$ are independent of the function $W(t)$, it follows that
\begin{equation} 
\PP(W(t_i) \geq 0, ~~ \forall 1 \leq i \leq N) = \EE \left [\frac{1}{N+1} \right ] = \sum_{k=0}^\infty \frac{1}{k+1} \frac{\alpha^k e^{-\alpha}}{k!} =
$$
$$
e^{-\alpha} \sum_{k=0}^\infty \frac{\alpha^k}{(k+1)!} = \frac{e^{-\alpha}}{\alpha} \sum_{k=1}^\infty \frac{\alpha^k}{k!} =
\frac{1 - e^{- \alpha}}{\alpha}.
\end{equation}
(recall that $N$ was a Poisson random variable with expectation $\alpha$).  \\ \\
We conclude the calculations as a lemma:
\begin{lemma} \label{lemsec2}
Let $B(\cdot)$ be a standard one-dimensional Brownian motion, let $W(\cdot)$ be a standard Brownian bridge on $[0,r]$ and let
$t_1,...,t_N$ be a Poisson point process on $[0,r]$ with intensity $\alpha$, all processes being independent. We have
\begin{equation} \label{walkest}
\P(B(t_i) \geq 0, ~~ \forall 1 \leq i \leq N) = \Psi(r \alpha) = \frac{1}{\sqrt{\pi r \alpha}} \Phi(r \alpha).
\end{equation}
where
\begin{equation} \label{defpsi}
\Psi(t) = \frac{2}{\pi \sqrt{t}} \int_0^{\sqrt{t}} \frac{e^{-x^2} dx}{\sqrt{1-x^2/t}}, ~~ \Phi(t) = \frac{2}{\sqrt{\pi}} \int_0^{\sqrt{t}} \frac{e^{-x^2} dx}{\sqrt{1-x^2/t}}.
\end{equation}
Moreover,
\begin{equation} \label{bridgeest}
\PP(W(t_i) \geq 0, ~~ \forall 1 \leq i \leq N) = \frac{1 - e^{- r \alpha}}{r \alpha}.
\end{equation}
\end{lemma}
\bigskip
Next, we will need the following estimate:
\begin{lemma} \label{exppos}
Let $B(\cdot)$ be a standard one-dimensional Brownian motion and let $T=(t_1,...,t_N)$ be a Poisson point process with intensity $\alpha$ on the interval $[0,L]$, independent from the Brownian motion.Define
$$
A = \left \{B(t_i) \geq 0, ~ \forall 1 \leq i \leq N \right \}.
$$
We have
$$
\EE \left [ \mathbf{1}_{  B(L) > 0 } B(L) \mathbf{1}_{A}  \right ] =   \frac{1}{\sqrt{2 \alpha}} \mathbf{erf} \left (\sqrt{\alpha L} \right ) + \frac{e^{-\alpha L} - 1 }{\sqrt{2 \pi L} \alpha}
$$
where $\mathbf{erf}$ is the error function, defined as
$$
\mathbf{erf}(x) = \frac{2}{\sqrt \pi} \int_0^x e^{-t^2} dt.
$$
\end{lemma}
\bigskip
\emph{Proof:} \\
The proof follows the same lines as in \cite[Page 215]{MP}. Let $u(t,x)$ be the function satisfying
\begin{equation} \label{heat}
u_t(t,x) = \frac{1}{2} u_{xx}(t,x) - U(x) u(t,x), ~~u(0,x) = x \mathbf{1}_{\{x \geq 0\}},
\end{equation}
where $U(x) = \alpha \mathbf{1}_{\{x < 0 \}}$. By the Feynman-Kac formula, one has
\begin{equation} \label{FK}
u(t,0) = \EE \left [ \mathbf{1}_{  B(t) > 0  } B(t) \exp \left( - \int_0^t U(B(t)) dt  \right ) \right ].
\end{equation}
(see e.g., \cite[Theorem 7.43, page 214]{MP} for a proof of this formula in the special case that $u(0,x) \equiv 1$. A straightforward adaptation of this proof may extend the formula to the boundary condition in hand). By the definition of the Poisson process and by the independence of $B(t)$ and the Poisson process, one has almost surely
$$
\PP(A | \mathcal{F}_B ) = \exp \left ( - \int_0^t U(B(t)) dt \right ).
$$
where $\mathcal{F}_B$ denotes the $\sigma$-algebra generated by the Brownian motion $B(\cdot)$.
Consequently,
$$
\EE( \mathbf{1}_{  B(t) > 0  } B(t) \mathbf{1}_A | \mathcal{F}_B ) = \mathbf{1}_{  B(t) > 0  } B(t) \exp \left ( - \int_0^t U(B(t)) dt \right )
$$
almost surely. Taking expectation on both sides and using (\ref{FK}) yields
$$
u(t,0) = \EE \left [\mathbf{1}_{  B(t) > 0  } B(t) \mathbf{1}_{A}  \right ].
$$
Our goal is therefore to estimate $u(L,0)$. For $\rho > 0$, we define
$$
g(x) = \int_0^\infty e^{- \rho t} u(t,x) dt,
$$
the Laplace transform of $u(\cdot, x)$. Integration by parts yields
$$
\left . e^{- \rho t} u(t,x) \right |_0^\infty = - \rho \int_0^\infty e^{- \rho t} u(t,x) dt + \int_0^\infty e^{- \rho t} u_t(t,x) dt,
$$
so, using (\ref{heat}),
$$
- u(0,x) = - \rho g(x) + \frac{1}{2} g''(x) - U(x) g(x).
$$
In other words,
\begin{equation} \label{eqsg}
- x + \rho g(x) - \frac{1}{2} g''(x) = 0, ~~ \forall x > 0,
\end{equation}
$$
(\rho + \alpha) g(x) - \frac{1}{2} g''(x) = 0, ~~ \forall x < 0.
$$
Next, we claim that the function $g(x) / (1 + x^2)$ is bounded. Indeed, consider the solution to the equations
\begin{equation}
w_t(t,x) = \frac{1}{2} w_{xx}(t,x), ~~w(0,x) = x^2 + 1.
\end{equation}
Another application of the Feynman-Kac formula yields
$$
w(t,x) = \EE \left [ (B(t) + x)^2 + 1 \right ]
$$
and, likewise
$$
u(t,x) = \EE \left [ \mathbf{1}_{  (B(t) + x) > 0  } (B(t) + x) \exp \left( - \int_0^t U(B(t) + x) dt  \right ) \right ] \leq
$$
$$
\EE \left [ |B(t) + x| \right ].
$$
These two equations easily imply that $w(t,x) \geq u(t,x)$ for all $x \in \RR$ and $t>0$. Now, it is easy to check that $w(x,t) = x^2 + 1 + t$, therefore $g(x) \leq \int_0^\infty e^{- \rho t} w(x,t) dt = 1/\rho^2 + 1 + x^2$. It follows that $g(x)$ cannot grow exponentially with 
$x$. The only solution to the equations (\ref{eqsg}) which satisfies this is
$$
g(x) = \frac{x}{\rho} + A e^{- \sqrt{2 \rho} x}, ~~ \forall x > 0,
$$
$$
g(x) = B e^{ \sqrt{2 (\rho + \alpha)} x}, ~~ \forall x < 0,
$$
for some $A,B \in \RR$. The function $g(x)$ should be continuously differentiable at $0$, thus, by matching derivatives, we attain
$$
A=B; ~~ 1/\rho - \sqrt{2 \rho} A = \sqrt{2 (\rho + \alpha)} B
$$
which gives
\begin{equation} \label{laplace1}
g(0) = A = B = \frac{1}{\sqrt{2} \rho (\sqrt{\rho} + \sqrt{\rho + \alpha}) } = \frac{\sqrt{\rho + \alpha} - \sqrt{\rho} }{ \sqrt{2} \rho \alpha }.
\end{equation}
According to Lerch's theorem, the Laplace transform is unique in the sense that if a continuous function $F(t)$ satisfies
\begin{equation} \label{invlaplace}
\int_0^\infty e^{-\rho t} F(t) dt = \frac{\sqrt{\rho + \alpha} - \sqrt{\rho} }{ \sqrt{2} \rho \alpha }, ~~ \forall \rho > 0,
\end{equation}
then necessarily one has
\begin{equation} \label{conslaplace}
F(t) = u(t,0) = \EE \left [ B(t) \mathbf{1}_{A}  \right ].
\end{equation}
Our goal is therefore to find a function $F(t)$ solving equation (\ref{invlaplace}). To that end, fix $\gamma \in \RR$ and let $F(t)$ be a function which satisfies
$$
F'(t) = \gamma \left ( \frac{1}{t^{3/2}} e^{-\alpha t} - \frac{1}{t^{3/2}} \right ), ~ \forall t > 0
$$
and $F(0) = 0$. We have
$$
\int_0^\infty F(t) e^{- \rho t} dt =  \left . - \frac{1}{\rho} F(t) e^{- \rho t} \right |_0^\infty + \frac{1}{\rho} \int_0^\infty F'(t) e^{-\rho t} dt =
$$
(plugging in the definition of $F$)
$$
 \frac{\gamma}{\rho} \int_0^\infty \frac{1}{t^{3/2}} \left( e^{-\alpha t} - 1 \right ) e^{- \rho t } dt =
$$
(integrating by parts)
$$
- \frac{2 \gamma }{\rho} \left . \left ( \frac{1}{\sqrt t} \left( e^{-\alpha t} - 1 \right ) e^{- \rho t } \right  ) \right |_0^\infty +
\frac{2 \gamma }{\rho} \int_0^\infty \frac{1}{\sqrt t} \left(- (\alpha + \rho) e^{- (\alpha+ \rho) t} + \rho e^{- \rho t }  \right ) dt.
$$
Now, a simple calculation shows that for every $\delta > 0$, one has
$$
\int_0^\infty \frac{e^{-\delta t}}{\sqrt t} dt = \sqrt{\frac{\pi}{\delta}}.
$$
So, the two above equations together yield
$$
\int_0^\infty F(t) e^{- \rho t} dt = \frac{2 \gamma \sqrt{\pi} }{\rho} \left ( - \sqrt{\alpha + \rho} + \sqrt{\rho} \right ).
$$
Finally, choosing $\gamma = \frac{- 1}{2 \sqrt{2 \pi} \alpha}$, gives
$$
\int_0^\infty F(t) e^{- \rho t} dt = \frac{\sqrt{\rho + \alpha} - \sqrt{\rho} }{ \sqrt{2} \rho \alpha }
$$
which means that (\ref{invlaplace}) is satisfied, and according to (\ref{conslaplace}) we conclude that
$$
\EE \left [ B(L) \mathbf{1}_{A}  \right ] = \frac{1}{2 \sqrt{2 \pi} \alpha} \int_0^L \frac{1}{s^{3/2}} \left ( 1 - e^{-\alpha s} \right ) ds = 
$$
$$
\frac{1}{\sqrt{2 \pi} \alpha} \left ( \sqrt{\pi \alpha} \mathbf{erf} \left (\sqrt{\alpha L} \right ) + \frac{e^{-\alpha L} - 1 }{\sqrt{L}} \right ) = 
$$
$$
  \frac{1}{\sqrt{2 \alpha}} \mathbf{erf} \left (\sqrt{\alpha L} \right ) + \frac{e^{-\alpha L} - 1 }{\sqrt{2 \pi L} \alpha}
$$
The proof is complete. \qed 

\section{A formula for the facets}

The goal of this section is to derive a formula which will serve as a central ingredient in our theorems.  \\ \\
We begin with some notation. Let $\Delta_n$ be the $n$-dimensional simplex, namely
$$
\Delta_n = \left \{ (r_1,...,r_n) \in [0,1]^n ; ~ \sum_{i=1}^n r_i \leq 1 \right \}.
$$
For a point $r = (r_1,...,r_n) \in \Delta_n$, define
$$
s_i(r) = \sum_{j=1}^i r_j, ~~ \forall 1 \leq i \leq n
$$
and
\begin{equation} \label{defs}
s(r) = (s_1(r), ..., s_n(r)).
\end{equation}
Next, for $r \in \Delta_n$ we define
$$
F_r = \Conv \bigl (B(s_1(r)),...,B(s_n(r)) \bigr)
$$
which is almost surely an $(n-1)$-dimensional simplex. Let $n_r$ be a unit vector normal to $F_r$ chosen such that $\langle n_r, B(s_1(r)) \rangle \geq 0$, and write
$$
V(r) = \Vol_{n-1}(F_r), ~~ H(r) = \langle n_r, B(s_1(r)) \rangle.
$$

Next, we define two point processes on $\Delta_n$. For a Borel subset $A \subset \Delta_n$ we define
\begin{equation} \label{defq}
q(A) =  \# \{r \in A; ~ F_r \mbox{ is a facet in the boundary of } K \}.
\end{equation}
and
$$
q_\alpha(A) = \# \{r \in A; ~ F_r \mbox{ is a facet in the boundary of } K_\alpha  \}.
$$
We also need the definition of the point process
$$
w_\alpha(A) = \# \left \{r \in A; ~~ \eta(r) \subset \Lambda_\alpha \right  \}
$$
where
$$
\eta(r) = \bigcup_{i=1}^n \left \{ s_i(r) \right \}
$$
which we can think of as points $r \in \Delta_n$ which are candidates to be facets of $K_\alpha$ in the sense that all their vertices are in the random walk. Observe that, since $K_\alpha$ is a polytope, one has almost surely
\begin{equation} \label{qleqw}
q_\alpha(A) \leq w_\alpha(A), ~~ \forall \alpha > 0, ~ A \subset \Delta_n.
\end{equation}
Moreover, we define the (deterministic) measures
$$
\mu(\cdot) = \EE [q(\cdot)], ~ \mu_\alpha(\cdot) = \EE \left [q_\alpha(\cdot) \right ] \mbox {  and  } \nu_\alpha(\cdot) = \EE [w_\alpha(\cdot)].
$$
The observation (\ref{qleqw}) implies $\mu_\alpha \ll \nu_\alpha$. Therefore, we may denote $$p_\alpha(r) = \frac{d \mu_\alpha}{d \nu_\alpha}(r), ~~ \forall r \in \Delta_n.$$
Let us try to understand how to calculate $p_\alpha(r)$. To this end, we need a few more definitions.

Let $\mathcal{F}_B$ be the $\sigma$-algebra generated by the Brownian motion $B(\cdot)$ (so that a random variable is measurable with respect to $\mathcal{F}_B$ if and only if
it does not depend on the point process $\Lambda$). We denote by $\mathcal{P}$ the space of all finite subsets of $[0,1]$. Next, for all $r \in \Delta_n$, we define the corresponding Palm measure
$$
P_{r, \alpha} (\cdot) = \PP_\Lambda \left  ( \Lambda_\alpha \cup \eta(r) \in \cdot \right ).
$$

Now, let $g: \Delta_n \times \mathcal{P} \to \RR$ be a function such that for all $r \in \Delta_n$ and $\phi \in \mathcal{P}$, the random variable $g(r, \phi)$ is measurable with respect to $\mathcal{F}_B$. According to the reduced Campbell-Little-Mecke formula (see e.g. \cite{SW}, section 3), one has for all $\alpha > 0$,
\begin{equation}
\EE \left [ \int_{\Delta_n} g(r, \Lambda_\alpha) w_\alpha(dr)  \right ] = \EE_B \left [  \int_{\Delta_n} \int_{\mathcal{P}}  g(r, \phi) P_{r, \alpha}(d \phi) \nu_\alpha(dr)   \right ].
\end{equation}
For $r \in \Delta_n$, define $M_r \subset \mathcal{P}$ to be the random set, measurable with respect to $\mathcal{F}_B$, satisfying
$$
\Lambda_\alpha \in M_r \Leftrightarrow F_r \mbox{ is a facet in the boundary of } \partial K_\alpha.
$$
Let $f: \Delta_n \to \RR$ be a function such that for all $r \in \Delta_n$, $f(r)$ is a random variable which is measurable with respect to $\mathcal{F}_B$. Then by
taking $g(r,\phi) = f(r) \mathbf{1}_{\{\phi \in M_r\} }$ in the previous formula and using Fubini's theorem, we obtain
$$
\EE \left [  \int_{\Delta_n} f(r) d q_\alpha (r)  \right ] = \EE_B \left [  \int_{\Delta_n} \int_{\mathcal{P}}  f(r) \mathbf{1}_{ \{\phi \in M_r \} } P_{r, \alpha}(d \phi) \nu_\alpha(dr)   \right ] = 
$$
$$
\int_{\Delta_n} \EE_B \left [ f(r) \int_{\mathcal{P}}  \mathbf{1}_{ \{\phi \in M_r \} } P_{r, \alpha}(d \phi)  \right ] \nu_\alpha(dr).
$$
At this point, it is convenient to define the events
$$
E_\alpha(r) := \left \{ F_r \mbox { is a facet in the boundary of } Conv \left ( \left (  \bigcup_{i=1}^n B(s_i(r)) \right ) \cup K_\alpha \right  ) \right \}.
$$
By the definition of $P_{r,\alpha}$ and $M_r$, we have 
$$
\int_{\mathcal{P}} \mathbf{1}_{ \{\phi \in M_r \} } P_{r, \alpha}(d \phi) = \EE \left [ \mathbf{1}_{E_\alpha(r)} \mid \mathcal{F}_B  \right ]
$$
almost surely with respect to $B(\cdot)$. A combination the two above formulas finally gives
\begin{equation} \label{superfubini}
\EE \left [  \int_{\Delta_n} f(r) d q_\alpha (r)  \right ] = \int_{\Delta_n} \EE \left [f(r) \mathbf{1}_{E_\alpha(r)} \right ] \nu_\alpha(d r).
\end{equation}
By taking $f(r) = 1$ we get that
\begin{equation} \label{prdef}
p_\alpha(r) = \P(E_\alpha(r)).
\end{equation}
Next, we would like to understand the measure $\nu_\alpha$. To that end, let $s = (s_1,...,s_n)$ and $\epsilon > 0$ be such that $s_i - s_{i-1} > \epsilon$ for all $2 \leq i \leq n$. 
Define 
$$
Q = s^{-1} (\{(x_1,...,x_n); ~ x_i \in [s_i, s_i + \epsilon], ~\mbox{for }i=1,..,n  \})
$$
where $s^{-1} (\cdot)$ is the inverse of the function defined in (\ref{defs}). Then, by the independence of the number of Poisson points on disjoint intervals,
$$
\nu (Q) = \EE \left [\prod_{i=1}^n \#\{j;~t_j \in [s_i, s_i + \epsilon] \} \right ] = (\epsilon \alpha)^n.
$$
By the $\sigma$-additivity of $\nu$, it follows that for a measurable $A \subset Int(\Delta_n) $,
$$
\nu(A) = \alpha^n \Vol_n(s(A)) = \alpha^n \Vol_n(A).
$$
where in the last equality we use the fact that the Jacobian of the function $r \to s(r)$ is identically one. We learn that, in fact,
$d \nu_\alpha = \alpha^n dr$.
In view of this identity, equation (\ref{superfubini}) becomes
\begin{equation} \label{superfubini2}
\EE \left [  \int_{\Delta_n} f(r) d q_\alpha (r)  \right ] = \alpha^n \int_{\Delta_n} \EE \left [f(r) \mathbf{1}_{E_\alpha(r)} \right ] dr.
\end{equation}
The above formula will play a central role in our proofs. It will serve us to find the expectation of several quantities of interest. For instance, in order to calculate the volume or the surface area of $K$, we observe that
\begin{equation} \label{voleq1}
\Vol_n (K_\alpha) = \sum_{r \in \Delta_n} \mathbf{1}_{\{ F_r \mbox{ is a facet of } K_\alpha \} } \Vol_n(\Conv(\{0\}, F_r) = 
\end{equation}
$$
\frac{1}{n} \int_{\Delta_n} V(r) H(r) d q_\alpha(r)
$$
and
\begin{equation} \label{volsa1}
\Vol_{n-1} (\partial K_\alpha) = \int_{\Delta_n} V(r) d q_\alpha(r).
\end{equation}
Using equation (\ref{superfubini2}), we get
\begin{equation} \label{volcalcidea}
\EE[\Vol_n(K_\alpha)] = \frac{\alpha^n}{n} \int_{\Delta_n} \EE \left [V(r) H(r)  \mathbf{1}_{E_\alpha(r)} \right ] dr.
\end{equation}
In a similar way, we obtain the following formula for the surface area:
\begin{equation} \label{sacalc}
\EE[\Vol_{n-1}(\partial K_\alpha)] = \alpha^n \int_{\Delta_n} \EE \left [V(r) \mathbf{1}_{E_\alpha(r)} \right ] dr.
\end{equation}

Next, we would like to derive more explicit expressions for the expectations on the right hand side of the two last formulae. The next lemma follows lines analogous to the ones developed in (\cite{E1}):

\begin{lemma} \label{lemfacets}
For all $\alpha > 0$ and all $r \in \Delta_n$, one has
\begin{equation} \label{facets2}
\PP(E_\alpha(r)) = 
2 \left (\prod_{j=2}^{n} \frac{1 - e^{-\alpha r_j}}{\alpha r_j} \right )  \Psi(\alpha r_1) \Psi(\alpha r_{n+1}).
\end{equation}
where $\Psi$ is defined as in Lemma \ref{lemsec2} and $r_{n+1} = 1 - \sum_{i=1}^n r_i$. Moreover, the event $E_\alpha(r)$ is independent from the equivalence class of $F_r$ (up to translations), and one has
\begin{equation} \label{voleq0}
\EE \left [V(r) H(r) \mathbf{1}_{E_\alpha(r)} \right ] = \EE[V(r)] \EE_B \left [H(r) \mathbf{1}_{E_\alpha(r)} \right ] = 
\end{equation}
$$
2 \EE[V(r)] \left (\prod_{j=2}^{n} \frac{1 - e^{-\alpha r_j}}{\alpha r_j} \right )  \Psi(\alpha r_{n+1}) \left ( \frac{1}{\sqrt{2 \alpha}} \mathbf{erf} \left (\sqrt{\alpha r_1} \right ) + \frac{e^{-\alpha r_1} - 1 }{\sqrt{2 \pi r_1} \alpha} \right ).
$$
\end{lemma}
\emph{Proof:} \\
Our first goal is to write the event $E_r$ as the product of independent events whose probabilities will be calculated using the formulas derived in section 2. The idea which allows us to do this is the following: the representation theorem for the Brownian bridge suggests that we may equivalently construct $B(t)$ by first generating the differences $B(s_j)-B(s_{j-1})$ as independent Gaussian random vectors, and then "fill in" the gaps between them by generating a Brownian motion up to $B(s_1)$, a Brownian bridge for each $1 < j\leq n$, and a "final" Brownian motion between $B(s_n)$ and $B(1)$, all of the above independent from each other. To make it formal, fix $r \in \Delta_n$ and define $s = s(r)$. For all $i$, $2 \leq i \leq n$, we write
$$
D_i = B(s_i) - B(s_{i-1})
$$
and define $C_i:[s_{i-1}, s_i] \to \RR^n$ by
$$
C_i(t) = B(t) - B(s_{i-1}) - \frac{t - s_{i-1}}{s_i - s_{i-1}} (B(s_i) - B(s_{i-1})),
$$
the bridges that correspond to the intervals $[s_{i-1}, s_i]$. Finally, we define two functions $B_0:[0,s_1] \to \RR^n$ and
$B_f:[s_n,1] \to \RR^n$ by $B_0(t) = B(s_1 - t) - B(s_1)$ and $B_f(t) = B(t) - B(s_n)$. By the independence of the differences of a Brownian motion on
disjoint intervals and by the representation theorem for the Brownian bridge, it follows that the variables $\{D_i\}_{i=2}^n, \{C_i\}_{i=2}^n, B_0, B_f$ are all
independent, each $C_i$ being a Brownian bridge and $B_0$ and $B_f$ being Brownian motions. \\ \\
Define
$$
\tilde C_i = \langle C_i, n_s \rangle, ~~\forall 2 \leq i \leq n
$$
and also $\tilde B_0 = \langle B_0, n_s \rangle$ and $\tilde B_f = \langle B_f, n_s \rangle$. Since $n_s$ is fully determined
by $\{D_i\}_{i=2}^n$, it follows that $\{\tilde C_i\}_{i=2}^n$, $\tilde B_0$ and $\tilde B_f$ are independent. Observe that for all $2 \leq i \leq n$, $\tilde C_i$ is a one-dimensional Brownian bridge fixed to be zero at its endpoints, and $\tilde B_0$ and $\tilde B_f$ are one dimensional Brownian motions starting from the origin. \\ \\
A moment of reflection reveals that the event $E_\alpha(s)$ is reduced to the intersection of the following conditions, for each possible direction
of $n_s$ with respect to $F_s$, \\ \\
(i) For all $2 \leq i \leq n$, the function $\tilde C_i$ is non-positive at all points $t_j$ such that $s_{i-1} \leq t_j \leq s_{i}$. \\
(ii) The function $\tilde B_0$ is non-positive at all points $t_j$ such that $t_j < s_1$. \\
(iii) The function $\tilde B_f$ is non-positive at all points $t_j$ such that $s_n < t_j \leq 1$. \\ \\
As explained above, $\{\tilde C_i\}_{i=2}^n$, $\tilde B_0$ and $\tilde B_f$ are independent, thus we can estimate $p(r)$ using equations (\ref{walkest}) and (\ref{bridgeest}). We get
\begin{equation} \label{facets}
p_\alpha(r) = 2 \left (\prod_{j=2}^{n} \frac{1 - e^{-\alpha r_j}}{\alpha r_j} \right ) \Psi(\alpha r_1) \Psi(\alpha r_{n+1}).
\end{equation}
Note that the factor $2$ stems from that fact that $n_s$ has two possible directions. Formula (\ref{facets2}) is thus established. \\ \\
Next, we note that $H_r = \tilde B_0(0)$. Defining $\tilde E_\alpha(r)$ as the event that (iii) above holds, we use Lemma \ref{exppos} in order to learn that 
$$
\EE   \left [\mathbf{1}_{\tilde E_\alpha(r)} H_r  \right] = \frac{1}{\sqrt{2 \alpha}} \mathbf{erf} \left (\sqrt{\alpha r_1} \right ) + \frac{e^{-\alpha r_1} - 1 }{\sqrt{2 \pi r_1} \alpha}.
$$
Since (i) and (iii) above are independent from (ii) and from $\tilde B_0(0)$, we get
$$
\EE  \left [\mathbf{1}_{E_\alpha(r)} H(r)  \right ] =  
$$
$$
2 \left (\prod_{j=2}^{n} \frac{1 - e^{-\alpha r_j}}{\alpha r_j} \right ) \Psi(\alpha r_{n+1}) \left ( \frac{1}{\sqrt{2 \alpha}} \mathbf{erf} \left (\sqrt{\alpha r_1} \right ) + \frac{e^{-\alpha r_1} - 1 }{\sqrt{2 \pi r_1} \alpha} \right ).
$$
Finally, the equivalence class of the facet $F_s$ (up to translations) clearly only depends on the differences $D_i$ which are, as explained above, independent of $B_0, B_f$ and $C_i$. Thus, we learn that $V(r)$ is independent from the events (i)-(iii) above. In particular, 
$$
\EE[\mathbf{1}_{E_\alpha(r)} V(r) H_\alpha(r)] = \EE[V(r)] \EE[\mathbf{1}_{E_\alpha(r)} H_\alpha(r)].
$$
A combination of the last two equations gives (\ref{voleq0}).
 \qed \\ \\

\section{Expected volume and surface area}

The purpose of this section is to use the technique developed in the previous section in order to obtain a formula for the expected volume 
of $K$. The formula will be derived in the following way: first, we can find a formula for the expected volume of $K_\alpha$ by combining formula (\ref{volcalcidea}) and Lemma \ref{lemfacets}. Then, in order to find $\EE[\Vol_n(K)]$, we will establish the fact that the latter may be expressed as a limit of the former, by taking $\alpha \to \infty$. This fact is stated in a corollary below. \\ \\

We begin with some notation. For a convex body $L$ and for $\epsilon > 0$, we denote
$$
e(L,\epsilon) := \{x \in \RR^n | ~\exists y \in L, ~|y - x| \leq \epsilon \},
$$
the $\epsilon$-extension of $L$. For two convex bodies $L,T$, we denote
$$
d_H(L,T) := \inf \{\epsilon; ~ L \subset e(T,\epsilon) \mbox{ and } T \subset e(L,\epsilon) \},
$$
the Hausdorff distance between $L$ and $T$.

\begin{lemma}
Almost surely, one has
$$
\lim_{\alpha \to \infty} d_H(K_\alpha, K)  = 0.
$$
\end{lemma}
\emph{Proof:} \\
For $\theta \in \Sph$, we write $h_K(\theta) = \sup_{x \in K} \langle x, \theta \rangle$, the support function of $K$, and 
$H_\theta(t) = \{ x \in \RR^n;~ \langle x, \theta \rangle \leq t \}$. Note that
$$
K = \bigcap_{\theta \in \Sph} H_\theta(h_K(\theta)).
$$
Define
$$
K(\epsilon) = \bigcap_{\theta \in \Sph} H_\theta(h_K(\theta) - \epsilon).
$$
It is easy to verify that, $\lim_{\epsilon \to 0} d(K(\epsilon)), K) = 0$. Therefore, it is enough to show that almost surely, for every $\epsilon > 0$, there exist $\alpha_0$ such that for every $\alpha > \alpha_0$, one has 
$$
\PP(K_\alpha \supset K(\epsilon)) > 1 - \epsilon.
$$
To that end, for every $\theta \in \Sph$, define $K(\epsilon, \theta) = K \setminus H_\theta(h_k(\theta) - \epsilon)$, and
$$
D(\theta) = \left \{\theta' \in \Sph; ~ \langle x, \theta' \rangle > h_{K} - \epsilon, ~~ \forall x \in K(\epsilon / 2, \theta) \right \}.
$$
Evidently, $D(\theta)$ is an open set that contains $\theta$. Next, define
$$
r(\theta) = \sup \{r | ~ B(\theta,r) \subset D(\theta)\},
$$
where $B(\theta,r)$ is an open spherical cap of radius $r$, centered at $\theta$. The fact that $D(\theta)$ is open implies that $r(\theta) > 0$.
Moreover, one may verify that $r(\theta)$ is continuous with respect to $\theta$, and therefore attains a minimum, $r_0>0$. Now, take $\theta_1,...,\theta_M$ to be an $r_0$-net of the sphere. Suppose a set of
points $x_1,...,x_M \in K$ satisfy $x_i \in K(\epsilon / 2, \theta_i)$, and denote $C = \Conv(x_1,..,x_M)$. Then for all $\theta \in \Sph$, there
exists some $i$ such that $\theta \in B(\theta_i, r_0)$ which implies that $\langle x_i, \theta \rangle \geq h_K(\theta) - \epsilon$.
It follows that $h_C(\theta) \geq h_K(\theta) - \epsilon$, which implies that $K(\epsilon) \subset C$. It is therefore enough
to show that the following event has probability tending to 1:
$$
E = \bigcap_{1 \leq i \leq M} \{ \exists x \in K_{\alpha} \mbox{ such that } \langle x, \theta_i \rangle \geq h_{K}(\theta_i) - \epsilon / 2 \}.
$$
For all $1 \leq i \leq M$, define $T_i = B^{-1}(K(\epsilon / 2, \theta_i))$. By the continuity of $B$, this set has a positive measure,
which means that the probability that one of the points of the Poisson process is in $T_i$ tends to $1$ as $\alpha \to \infty$. By applying a union bound, it follows
that $\lim_{\alpha \to \infty} P(E) = 1$, and the lemma is proven. \qed \\ \\
As a direct corollary, we obtain
\begin{corollary} \label{corlimit}
Almost surely, one has
$$
\lim_{\alpha \to \infty} \Vol_n(K_\alpha) = \Vol_n(K),
$$
and
$$
\lim_{\alpha \to \infty} \Vol_{n-1}(\partial K_\alpha) = \Vol_{n-1}(\partial K).
$$
\end{corollary}

In view of the above corollary, the proof of Theorem \ref{thmvolume} is reduced to calculating $\lim_{\alpha \to \infty} \EE[\Vol_n(K_\alpha)]$. Recall formulae (\ref{volcalcidea}) and (\ref{voleq0}). 
The only ingredient we still need is $\EE[V(r)]$. The next lemma is a simple calculation.

\begin{lemma} \label{lemfacetvolume}
Let $r = (r_1,...,r_n) \in \Delta_n$. We have
$$
\EE[ V(r) ] = 2^{(n-1)/2} \frac{\Gamma((n+1)/2)}{\Gamma(n)} \prod_{i=1}^{n-1} \sqrt{r_{i+1}}.
$$
Furthermore, $V(r) \sim \prod_{i=1}^{n-1} \sqrt{r_{i+1}} X$ where $X$ is a random variable whose distribution does not depend on $r$.
\end{lemma}
\emph{Proof:}
Define $s = s(r)$ and
$$
v_i = B(s_{i+1}) - B(s_i)
$$
for $1 \leq i \leq n-1$. One has
$$
(n-1)! \Vol_{n-1} (F_r) = \left | \det \left (
\begin{matrix}
 v_1 \\
v_1 + v_2 \\
... \\
v_1 + ... + v_{n-1} \\
n_s
\end{matrix}
\right ) \right | = \left |
\det \left (
\begin{matrix}
 v_1 \\
v_2 \\
... \\
v_{n-1} \\
n_s
\end{matrix}
\right ) \right |.
$$
Let $\Gamma_1,..,\Gamma_{n-1}$ be independent standard Gaussian random vectors in $\RR^n$. By the independence of increments of the Brownian motion on disjoint intervals, we have 
$$
(v_1,..,v_{n-1}) \sim (\sqrt{r_2} \Gamma_1, ..., \sqrt{r_n} \Gamma_{n-1}).
$$
So,
$$
\Vol_{n-1} (F_r) \sim \frac{1}{(n-1)!} \prod_{i=1}^{n-1} \sqrt{r_{i+1}} \left | \det  \left (
\begin{matrix}
\Gamma_1 \\
... \\
\Gamma_{n-1} \\
n_r
\end{matrix}
\right ) \right |
$$
Denote $E_0 = \RR^n$ and $E_i = span \{\Gamma_1, \Gamma_2, ... \Gamma_i\}^{\perp}$. Thinking about the above determinant as the volume of the parallelepiped spanned by the vectors $\Gamma_i$ and using a simple induction to calculate this volume gives
$$
\left | \det  \left (
\begin{matrix}
\Gamma_1 \\
... \\
\Gamma_{n-1} \\
u
\end{matrix}
\right ) \right | = \prod_{i=1}^{n-1} \left | Proj_{E_{i-1}} \Gamma_i \right |.
$$
Observe that the dimension of $E_i$ is almost surely $n-i$. Let $u_1,...,u_n$ be vectors such that $\{ u_1,...,u_{n-i+1} \}$ is an orthonormal basis of $E_{i-1}$. One has
$$
\EE \left [ \left | Proj_{E_{i-1}} \Gamma_i  \right | \right ] = \EE \left [ \sqrt{\sum_{j=1}^{n-i+1} \langle \Gamma_i, u_j \rangle^2} \right ].
$$
Note that in the last equality we use the fact that $\Gamma_i$ is independent of $E_{i-1}$ and therefore the vectors $u_1,...,u_{n-i+1}$ can
be assumed constant. The above is just the first moment of the $\chi$-distribution with $(n-i+1)$ degrees of freedom, which is equal to
$$
\EE \left [ \left | Proj_{E_{i-1}} \Gamma_i  \right | \right ] = \sqrt{2} \frac{\Gamma((n-i+2)/2)}{\Gamma((n-i+1)/2)}.
$$
Since $\Gamma_i$ is independent of $E_{i-1}$, it also follows that the variables $|Proj_{E_{i-1}} \Gamma_i|$ are independent. Consequently,
$$
\EE \left [  \left | \det  \left (
\begin{matrix}
\Gamma_1 \\
... \\
\Gamma_{n-1} \\
u
\end{matrix}
\right ) \right |  \right ] = \prod_{i=1}^{n-1} \sqrt{2} \frac{\Gamma((n-i+2) / 2)}{\Gamma((n-i+1)/2)} = 2^{(n-1)/2} \Gamma((n+1)/2).
$$
We conclude,
$$
\EE[\Vol_{n-1} (F_r)] = 2^{(n-1)/2} \frac{\Gamma((n+1)/2)}{\Gamma(n)} \prod_{i=1}^{n-1} \sqrt{r_{i+1}}.
$$
\qed \\ \\

We now have all the ingredients we need. Plugging the result of the above lemma in (\ref{voleq0}), we attain
\begin{equation} 
\EE \left  [ \mathbf{1}_{E_\alpha(r)} H(r) V(r) \right ] = 
\end{equation}
$$
2 \left (\prod_{j=2}^{n} \frac{1 - e^{-\alpha r_j }}{\alpha r_j} \right ) \Psi(\alpha r_{n+1}) \EE[\mathbf{1}_{\tilde E_r} H_r] \EE[V_r] = 
$$
\begin{equation} 
\frac{\xi_n}{\alpha^n} \left (\prod_{j=2}^{n+1} \frac{1}{\sqrt{r_j}} \right ) \prod_{j=2}^{n} \left (1 - e^{-\alpha r_j} \right ) \left (\mathbf{erf} \left (\sqrt{\alpha r_1} \right ) + \frac{e^{-\alpha r_1} - 1 }{\sqrt{\pi r_1 \alpha}} \right ) \Phi(\alpha r_{n+1}).
\end{equation}
where we set $\xi_n = \frac{1}{\sqrt{\pi}} 2^{n/2} \frac{\Gamma((n+1)/2)}{\Gamma(n)}$. \\ \\
Plugging this into equation (\ref{volcalcidea}), we finally get
\begin{equation} \label{volcalc2}
\EE[\Vol_n(K_\alpha)] = 
\end{equation}
$$
\frac{\xi_n}{n} \int_{\Delta_n} \left (\prod_{j=2}^{n+1} \frac{1}{\sqrt{r_j}} \right ) \left ( \prod_{j=2}^{n} \left (1 - e^{-\alpha r_j} \right ) \right ) \Phi(\alpha r_{n+1}) \left (\mathbf{erf} \left (\sqrt{\alpha r_1} \right ) + \frac{e^{-\alpha r_1} - 1 }{\sqrt{\pi r_1 \alpha}} \right ) dr.
$$
As mentioned above, we aim at calculating $\EE[\Vol_n(K)]$ by means of taking $\alpha \to \infty$. To do this, we would like to use the dominated convergence theorem in order to take the limit inside the integral, so we end up with a simpler integrand. Let us inspect each term separately.
First, note that
$$
\prod_{j=2}^{n} \left (1 - e^{-\alpha r_j} \right ) \leq 1,
$$
and
$$
0 \leq \mathbf{erf} \left (\sqrt{\alpha r_1} \right ) + \frac{e^{-\alpha r_1} - 1}{\sqrt{\pi r_1 \alpha}} \leq \mathbf{erf} \left (\sqrt{\alpha r_1} \right ) \leq 1.
$$
(where the first inequality follows from Lemma \ref{exppos}).
Also, we estimate
\begin{equation} \label{estphi}
\Phi(t) = \frac{2}{\sqrt{\pi}} \int_0^{\sqrt{t}} \frac{e^{-x^2} dx}{\sqrt{1-x^2/t}} \leq
\end{equation}
$$
\frac{2}{\sqrt{\pi}} \left ( \int_0^{\sqrt{\frac{3}{4}t}} \frac{e^{-x^2} dx}{\sqrt{1/4}} dx + \int_{\sqrt{\frac{3}{4}t}}^{\sqrt{t}} \frac{e^{-\frac{3}{4} t}dx}{ {\sqrt{1-x^2/t}} } \right ) \leq 2 + e^{-\frac{3}{4} t} \sqrt{t} \sqrt{\pi} \leq 3,
$$
for all $t \geq 0$. We learn that the integrand in (\ref{volcalc2}) is positive and smaller than the term $3 \prod_{j=2}^{n+1} \frac{1}{\sqrt{r_j}}$, whose integral clearly converges, therefore we may use the dominated convergence theorem. For all $r \in Int(\Delta_n)$, we calculate
$$
\lim_{\alpha \to \infty} \left (\mathbf{erf} \left (\sqrt{\alpha r_1} \right ) + \frac{e^{-\alpha r_1} - 1}{\sqrt{\pi r_1 \alpha}} \right ) =
$$
$$
\lim_{\alpha \to \infty} \mathbf{erf} \left (\sqrt{\alpha r_1} \right ) = 1,
$$
and
$$
\lim_{\alpha \to \infty} \Phi(\alpha r_{n+1}) =  \lim_{t \to \infty} \frac{2}{\sqrt{\pi}} \int_0^{\sqrt{t}} \frac{e^{-x^2} dx}{\sqrt{1-x^2/t}} =
$$
$$
\frac{2}{\sqrt{\pi}} \int_0^\infty e^{-x^2} dx = 1.
$$
We attain
$$
\lim_{\alpha \to \infty} \EE[\Vol_n(K_\alpha)] = \frac{\xi_n}{n} \int_{\Delta_n} \prod_{j=2}^{n+1} \frac{1}{\sqrt{r_j}} dr =
$$
(interchanging between $r_1$ and $r_{n+1}$) 
$$
\frac{\xi_n}{n} \int_{\left \{ \sum_{i=1}^n r_i \leq 1 \right \} } \prod_{i=1}^{n} \frac{1}{\sqrt{r_i}} d r = 
$$
(substituting $t_i = \sqrt{r_i}$)
$$
\frac{\xi_n}{n} 2^n \int_{ \left \{ \sum_{i=1}^n t_i^2 \leq 1 \atop t_i \geq 0 \right \}} dt = \frac{\xi_n}{n} \int_{ \left \{ \sum_{i=1}^n t_i^2 \leq 1 \right \}} dt = \frac{\xi_n}{n} \Vol_n(\BBN) 
$$
where $\BBN := \{x \in \RR^n, |x| \leq 1 \}$, the unit Euclidean ball. \\ \\
The following formula is well-known:
$$
\Vol_n(\BBN) = \frac{\pi^{n/2}}{\Gamma \left (\frac{n}{2}+1 \right )}.
$$
We finally get
$$
\EE[\Vol_n(K)] = \frac{1}{\sqrt{\pi} n} (2 \pi)^{n/2} \frac{\Gamma((n+1)/2)}{\Gamma(n) \Gamma \left (\frac{n}{2} + 1 \right )} = \left (\frac{\pi}{2} \right )^{n/2} \frac{1}{\Gamma \left (\frac{n}{2} + 1 \right )^2}
$$
The computation of the surface area is completely analogous. By combining (\ref{sacalc}), (\ref{facets}) and Lemma \ref{lemfacetvolume} we get
$$
\EE[\Vol_{n-1}(\partial K)] = \int_{\Delta_n} \EE[V_r] \PP[1_{E_\alpha(r)} ] = 
$$
$$
\frac{2 }{\pi} 2^{(n-1)/2} \frac{\Gamma((n+1)/2)}{\Gamma(n)} \int_{\Delta_n} \left (\prod_{j=1}^{n+1} \frac{1}{\sqrt{r_j}} \right ) dr =
$$
(substituting $t_i = \sqrt{r_i}$)
$$
\frac{2}{\pi} 2^{(n-1)/2} \frac{\Gamma((n+1)/2)}{\Gamma(n)} 2^n \int_{\left \{ \sum_{i=1}^n t_i^2 \leq 1 \atop t_i \geq 0 \right \}} \frac{1}{\sqrt{1 - |t|^2}} dt =
$$
$$
\frac{2}{\pi} 2^{(n-1)/2} \frac{\Gamma((n+1)/2)}{\Gamma(n)} \int_{\BBN} \frac{1}{\sqrt{1 - |t|^2}} dt =
$$
$$
\Vol_{n-1}(\partial \BBN) \frac{2}{\pi} 2^{(n-1)/2} \frac{\Gamma((n+1)/2)}{\Gamma(n)} \int_0^1 \frac{x^{n-1}}{\sqrt{1 - x^2}} dx =
$$
$$
\frac{4 \pi^{n/2}}{\Gamma(n/2)} \frac{1}{\pi} 2^{(n-1)/2} \frac{\Gamma((n+1)/2)}{\Gamma(n)} \frac{\sqrt{\pi} \Gamma(n/2) }{2 \Gamma((n+1)/2)}  = \frac{2 (2 \pi)^{(n-1) / 2}}{\Gamma(n)}.
$$
We have established Theorem \ref{thmvolume}. \\ \\

The following proof of Corollary \ref{cormixed} was communicated to us by Christoph Th\"{a}le. \\
\emph{Proof of Corollary \ref{cormixed}:}

For a linear subspace $L \subset \RR^n$, define by $K|L$ the projection of $K$ onto $L$, in other words,
$$
K|L = \left \{x \in L; ~ \exists y \in L^\perp, ~x+y \in K \right \}.
$$
For all $1 \leq j \leq n$, let $E_j$ be some fixed $j$-dimensional subspace. According to Kubota's formula (see \cite[p. 222]{SW}), we have
$$
V_j(K) = c(n,j) \int_{ \textbf{SO}(n) } \Vol_j ( U (K) | E_j) m(dU)
$$
where $\textbf{SO}(n)$ denotes the special orthogonal group on $n$-dimensions, $m(\cdot)$ denotes the normalized Haar measure on this group and
$$
c(n,j) :=  \left ( n \atop j \right ) \frac{\Gamma \left (\frac{j}{2} + 1 \right ) \Gamma \left ( \frac{n-j}{2} + 1 \right ) }{ \Gamma(n/2 + 1)}.
$$
After taking expectation on both sides and using Fubini's theorem, the last formula becomes
$$
\EE V_j(K) = c(n,j) \int_{ \textbf{SO}(n) } \EE \left [\Vol_j ( U (K) | E_j) \right ] m(dU).
$$
Now, by the rotational invariance of the Brownian motion, for any fixed $U \in \mathbf{SO}(n)$, the body $K$ has the same distribution as $U(K)$, which tells us that
$$
\EE V_j(K) = c(n,j) \EE \left [\Vol_j ( K | E_j) \right ].
$$
Next, we make the following observation: since the coordinates of a Brownian motion, under any orthogonal basis, are independent and since the operation of taking the convex hull commutes with the operation of projecting a set onto a linear subspace, the random body $K|E_j$ has the same distribution as that of the convex hull of a $j$-dimensional Brownian motion embedded in $E_j$. According to Theorem \ref{thmvolume}, we therefore know that
$$
\EE \left [\Vol_j ( K | E_j) \right ] = \left (\frac{\pi}{2} \right )^{j/2} \frac{1}{\Gamma \left (\frac{j}{2} + 1 \right )^2}.
$$
A combination of the two last formulas yields the result of the corollary.
\qed

\section{The approximating polytope}

The goal of this section is to prove Theorem \ref{thmapprox}. Since in the previous section, we already calculated $\EE[\Vol_n(K)]$, the derivation
of the bounds in the theorem is reduced to obtaining estimates on $\EE[\Vol_n(K_\alpha)]$, which in turn, are reduced to obtaining respective
estimates on (\ref{volcalc2}). \\ \\
We begin with the lower bound.  Inspect equation (\ref{volcalc2}). Our first goal will be to show that the term involving the expression
$ \frac{e^{-\alpha r_1} - 1 }{\sqrt{\pi r_1 \alpha}}$ is small. We calculate, using (\ref{estphi}),
$$
\frac{\xi_n}{n} \int_{\Delta_n}  \left (\prod_{j=2}^{n+1} \frac{1}{\sqrt{r_j}} \right ) \left ( \prod_{j=2}^{n}  \left (1 - e^{-\alpha r_j} \right ) \right ) \Phi(\alpha r_{n+1}) \frac{1 - e^{-\alpha r_1}}{\sqrt{\pi r_1 \alpha}} dr \leq 
$$
$$
\frac{3 \xi_n}{n \sqrt{\alpha}} \int_{\Delta_n} \prod_{j=1}^{n+1} \frac{1}{\sqrt{r_j}} dr = \frac{3 \xi_n}{n \sqrt{\alpha}} \int_{\BBN} \frac{1}{\sqrt{1 - |x|^2}} dx
$$
(here, we used the substitution $r_i = x_i^2$). It is easy to verify that,
$$
\int_{\BBN} \frac{1}{\sqrt{1 - |x|^2}} dx \leq 4 \sqrt{n} \Vol_{n} (\BBN).
$$
So,
\begin{equation} \label{firstterm}
\frac{\xi_n}{n} \int_{\Delta_n}  \left (\prod_{j=2}^{n+1} \frac{1}{\sqrt{r_j}} \right ) \left ( \prod_{j=2}^{n}  \left (1 - e^{-\alpha r_j} \right ) \right ) \Phi(\alpha r_{n+1}) \frac{1 - e^{-\alpha r_1}}{\sqrt{\pi r_1 \alpha}} dr \leq  
\end{equation}
$$
\frac{12 \sqrt{n} }{\sqrt{\alpha}} \frac{\xi_n \Vol_n(\BBN)}{n}.
$$
Our next task is to estimate the remaining term. To that end, we note that $\Phi(x) \geq \mathbf{erf}(\sqrt x)$ and use the following well-known estimate:
$$
\mathbf{erf} (x) \geq 1 - e^{-x^2}.
$$
This estimate yields
$$
\frac{\xi_n}{n} \int_{\Delta_n} \left (\prod_{j=2}^{n+1} \frac{1}{\sqrt{r_j}} \right ) \left ( \prod_{j=2}^{n} \left (1 - e^{-\alpha r_j} \right ) \right ) \Phi(\alpha r_{n+1}) \mathbf{erf} \left (\sqrt{\alpha r_1} \right ) dr \geq
$$
$$
\frac{\xi_n}{n} \int_{\Delta_n} \left (\prod_{j=2}^{n+1} \frac{1}{\sqrt{r_j}} \right ) \left ( \prod_{j=1}^{n+1} \left (1 - e^{-\alpha r_j} \right ) \right ) dr = 
$$
When interchanging the roles of $r_1$ and $r_{n+1}$ in the above integral, the domain of integration $\Delta_n$ becomes $\{\sum_{i=1}^n r_i \leq 1\}$, and by using the substitution $t_i^2 = r_i$ for all $1 \leq i \leq n$, the domain of integration will become the unit ball, $\BBN$. In other words,
\begin{equation} \label{errfest}
\frac{\xi_n}{n} \int_{\Delta_n} \left (\prod_{j=2}^{n+1} \frac{1}{\sqrt{r_j}} \right ) \left ( \prod_{j=2}^{n} \left (1 - e^{-\alpha r_j} \right ) \right ) \Phi(\alpha r_{n+1}) \mathbf{erf} \left (\sqrt{\alpha r_1} \right ) dr \geq
\end{equation}
$$
\frac{\xi_n}{n} \int_{\sum_{i=1}^n r_i \leq 1} \left (\prod_{j=1}^{n} \frac{1}{\sqrt{r_j}} \right ) \left ( \prod_{j=1}^{n} \left (1 - e^{-\alpha r_j} \right ) \right ) \left ( 1 - e^{- \left (1 - \sum_{i=1}^n r_i \right ) \alpha }  \right ) dr = 
$$
$$
\frac{\xi_n}{n} \int_{\BBN} \left ( \prod_{j=1}^{n} \left (1 - e^{-\alpha t_j^2} \right ) \right ) \left ( 1 - e^{- \left (1 - |t|^2 \right ) \alpha }  \right ) dt.
$$
Define
$$
F(t) = \frac{\Vol_{n} \left ( \left \{ \left (1 - |x|^2 \right ) \leq t \right \} \cap \BBN \right  )}{\Vol_n(\BBN)}.
$$
A straightforward calculation gives
\begin{equation} \label{estF}
F \left ( \frac{t}{10 n} \right ) \leq t , ~~ \forall t \geq 0.
\end{equation}
So we can estimate
\begin{equation}
\int_{\BBN} e^{ \left (|x|^2 - 1 \right ) \alpha} dx \leq \int_{\BBN}  \frac{1}{1 + \left ( 1 - |x|^2 \right ) \alpha} dx =
\end{equation}
$$
\Vol_n(\BBN) \int_0^1 \frac{F'(t)}{1 + \alpha t} dt = \Vol_n(\BBN) \left ( \frac{F(1)}{1 + \alpha} + \alpha \int_0^1 \frac{F(t)}{(1 + \alpha t)^2} dt \right ).
$$
Using (\ref{estF}), we attain
$$
\int_{\BBN} e^{ \left (|x|^2 - 1 \right ) \alpha} dx \leq \Vol_n(\BBN) \left ( \frac{1}{\alpha} + 10 n \alpha \int_0^1 \frac{t}{(1 + \alpha t)^2} dt \right ) =
$$
$$
\Vol_n(\BBN) \left ( \frac{1}{\alpha} + \frac{10 n}{\alpha} \int_0^{\alpha} \frac{t}{(1 + t)^2} dt \right ) \leq \frac{n \Vol_n(\BBN)}{\sqrt \alpha},
$$
where in the last inequality, we use the legitimate assumption that $\alpha$ is greater than some universal constant. Combining this with (\ref{errfest}) yields
\begin{equation} \label{eq111}
\frac{\xi_n}{n} \int_{\Delta_n} \left (\prod_{j=2}^{n+1} \frac{1}{\sqrt{r_j}} \right ) \left ( \prod_{j=2}^{n} \left (1 - e^{-\alpha r_j} \right ) \right ) \Phi(\alpha r_{n+1}) \mathbf{erf} \left (\sqrt{\alpha r_1} \right ) dr \geq
\end{equation}
$$
\frac{\xi_n}{n} \left ( \int_{\BBN} \left ( \prod_{j=1}^{n} \left (1 - e^{-\alpha t_j^2} \right ) \right ) dt - \frac{n \Vol_n(\BBN)}{\sqrt \alpha} \right ).
$$
To estimate the right hand side, we will need the following inequality,
\begin{equation} \label{newineq}
\int_{\BBN} \left ( \prod_{j=1}^{n} \left (1 - e^{-\alpha t_j^2} \right ) \right ) dt \geq
\end{equation}
$$
\Vol_n(\BBN) \frac{\int_{\BBN} \left ( \prod_{j=1}^{n} \left (1 - e^{-\alpha t_j^2} \right ) \right ) e^{- 2n \sum_{i=1}^n t_i^2} dt}{\int_{\BBN} e^{- 2n \sum_{i=1}^n t_i^2} dt}.
$$
In order to see why this inequality holds, we first note that the term $\prod_{j=1}^{n} \left (1 - e^{-\alpha t_j^2} \right )$ is monotone on rays of the form 
$\{s (t_1,...,t_n); ~ s \geq 0 \}$ and we recall the simple fact that for any two integrable monotone functions $f,g: [0,1] \to \RR^+$ and for any probability measure $\mu$ on $[0,1]$ one has
$$
\int_0^1 f(s) g(s) d \mu(s) \geq \int_0^1 f(s) d \mu(s) \int_0^1 g(s) d \mu(s).
$$
Fix $\theta = (\theta_1,...,\theta_n) \in \partial \BBN$. By taking $\mu$ to be the measure satisfying $\frac{d \mu}{ds} =n s^{n-1}$ and defining
$$
f(s) = e^{- 2n s \sum_{i=1}^n \theta_i^2}, ~~ g(s) = \prod_{j=1}^{n} \left (1 - e^{-\alpha s \theta_j^2} \right )
$$
we get
$$
\frac{\int_0^1 s^{n-1} \left ( \prod_{j=1}^{n} \left (1 - e^{-\alpha s \theta_j^2} \right ) \right ) e^{- 2n \sum_{i=1}^n s \theta_i^2} ds}{\int_0^1 s^{n-1} e^{- 2n \sum_{i=1}^n \theta_i^2} ds} \leq
$$
$$
n \int_0^1 s^{n-1} \left ( \prod_{j=1}^{n} \left (1 - e^{-\alpha s \theta_j^2} \right ) \right ) ds.
$$
Now, by rotational symmetry,
$$
\Vol_{n-1}(\partial \BBN) \int_0^1 s^{n-1} e^{- 2n \sum_{i=1}^n \theta_i^2} ds = \int_{\BBN} e^{- 2n \sum_{i=1}^n t_i^2} dt
$$
and therefore 
$$
\frac{\int_0^1 s^{n-1} \left ( \prod_{j=1}^{n} \left (1 - e^{-\alpha s \theta_j^2} \right ) \right ) e^{- 2n \sum_{i=1}^n s \theta_i^2} ds}{\int_{\BBN} e^{- 2n \sum_{i=1}^n t_i^2} dt} \leq
$$
$$
\frac{n}{ \Vol_{n-1}(\partial \BBN)} \int_0^1 s^{n-1} \left ( \prod_{j=1}^{n} \left (1 - e^{-\alpha s \theta_j^2} \right ) \right ) ds = 
$$
$$
\frac{1}{\Vol_{n}(\BBN)} \int_0^1 s^{n-1} \left ( \prod_{j=1}^{n} \left (1 - e^{-\alpha s \theta_j^2} \right ) \right ) ds.
$$
Integration of both sides of the last equation with respect to $\theta$ on $\partial \BBN$ finally proves \eqref{newineq}.

Next, since the term $e^{- 2n \sum_{i=1}^n t_i^2}$ is proportional to the density of a standard Gaussian random vector $\Gamma=(\Gamma_1,...,\Gamma_n)$, we get
$$
\int_{\BBN} \left ( \prod_{j=1}^{n} \left (1 - e^{-\alpha t_j^2} \right ) \right ) dt \geq
$$
$$
\Vol_n(\BBN) \left . \EE \left [\prod_{j=1}^{n} \left (1 - e^{- \frac{\alpha}{4n} \Gamma_j^2} \right ) ~~ \right | |\Gamma|^2 \leq 4n \right ].
$$
A calculation gives $\PP(|\Gamma|^2 \leq 4n) > 1- e^{-n}$. So we get
$$
\int_{\BBN} \left ( \prod_{j=1}^{n} \left (1 - e^{- \alpha t_j^2} \right ) \right ) dt \geq \Vol_n(\BBN) \left (\prod_{j=1}^n \EE \left [1 - e^{-\frac{\alpha}{4n} \Gamma_j^2} \right ] - e^{-n} \right )= 
$$
$$
\Vol_n(\BBN) \left ( \left ( 1 - \frac{1}{\sqrt{\frac{\alpha}{2n} + 1}} \right )^n - e^{-n} \right ) \geq \Vol_n(\BBN) \left (1 - 2 \frac{n^{3/2}}{\sqrt{\alpha}} - e^{-n} \right )
$$
The above estimate combined with (\ref{eq111}) gives
$$
\frac{\xi_n}{n} \int_{\Delta_n} \left (\prod_{j=2}^{n+1} \frac{1}{\sqrt{r_j}} \right ) \left ( \prod_{j=2}^{n} \left (1 - e^{-\alpha r_j} \right ) \right ) \Phi(\alpha r_{n+1}) \mathbf{erf} \left (\sqrt{\alpha r_1} \right ) dr \geq
$$
$$
\frac{\xi_n \Vol_n(\BBN) }{n} \left (1 - 2 \frac{n^{3/2}}{\sqrt{\alpha}} - e^{-n} - \frac{n}{\sqrt \alpha} \right ) \geq 
\frac{\xi_n \Vol_n(\BBN) }{n} \left (1 - 3 \frac{n^{3/2}}{\sqrt{\alpha}} - e^{-n} \right ).
$$
Finally, the last equation along with (\ref{volcalc2}) and (\ref{firstterm}) give
$$
\EE[\Vol_n(K_\alpha)] \geq \frac{\xi_n \Vol_n(\BBN) }{n} \left (1 - 15 \frac{n^{3/2}}{\sqrt{\alpha}} - e^{-n} \right )
$$
Since $K_\alpha \subset K$ almost surely, one has $\EE[\Vol_n(K \setminus K_\alpha)] = \EE[\Vol_n(K)] - \EE[\Vol_n(K_\alpha)]$. Therefore, combining the above bound
with Theorem \ref{thmvolume} establishes (\ref{aproxlower}). \\ \\

We continue with proving the bound (\ref{aproxupper}). Formula (\ref{volcalc2}) and the bound (\ref{estphi}) suggest that
\begin{equation} \label{lowereq1}
\EE[\Vol_n(K_\alpha)] \leq \frac{3 \xi_n}{n} \int_{\Delta_n} \left (\prod_{j=2}^{n+1} \frac{1}{\sqrt{r_j}} \right ) \left ( \prod_{j=2}^{n} \left (1 - e^{-\alpha r_j} \right ) \right ) dr = 
\end{equation}
(interchanging between $r_1$ and $r_{n+1}$ and substituting $t_i^2 = r_i$)
$$
\frac{3 \xi_n}{n} \int_{\BBN} \left ( \prod_{j=2}^{n} \left (1 - e^{-\alpha t_j^2} \right ) \right ) dt.
$$
Let $X = (X_1,...,X_n)$ be uniformly distributed in $\BBN$. It is straightforward to verify that
$$
\PP \left (|X_1| \leq \frac{t}{n \sqrt n} \right ) \geq \frac{t}{4n}, ~~ \forall 0 \leq t \leq \frac{n}{10}.
$$
We observe that for all $k>1$ and for all $0 \leq t \leq n/10$, one has
$$
\PP \left . \left (|X_k| \leq  \frac{t}{n \sqrt n} ~\right | ~ |X_1| >  \frac{t}{n \sqrt n}, |X_2| >  \frac{t}{n \sqrt n}, ..., |X_{k-1}| >  \frac{t}{n \sqrt n} \right ) \geq
$$
$$
 \PP \left (|X_k| \leq  \frac{t}{n \sqrt n} \right ).
$$
It follows that
$$
\PP \left (\min_k |X_k| > \frac{t}{n \sqrt n} \right ) = 
$$
$$
P(|X_1| > \frac{t}{n \sqrt n}) \prod_{k=2}^n \PP \left (|X_k| \leq t ~ | ~ |X_1| > \frac{t}{n \sqrt n}, ..., |X_{k-1}| > \frac{t}{n \sqrt n} \right  )  \leq e^{-t/4},
$$
for all $0 \leq t \leq n/10$. The last inequality, combined with (\ref{lowereq1}), implies that
$$
\EE[\Vol_n(K_\alpha)] \leq 
$$
$$
\frac{3 \xi_n \Vol_n(\BBN)}{n} \left ( e^{-t/4} + (1 - e^{-t/4}) \left (1 - e^{- \alpha t^2 / n^3 } \right ) \right ) \leq
$$
$$
\frac{3 \xi_n \Vol_n(\BBN)}{n} \left (1 + e^{-t/4} - e^{- \alpha t^2 / n^3 } \right ) \leq \frac{3 \xi_n \Vol_n(\BBN)}{n} \left (e^{-t/4} + \alpha t^2 / n^3 \right ),
$$
for all $0 \leq t \leq n/10$. Using the assumption $\alpha < \frac{n^3}{8}$ and choosing $t = 4 \log (n^3 / \alpha)$ gives
$$
\EE[\Vol_n(K_\alpha)] \leq 100 \frac{\alpha}{n^3} \log^2 \left (\frac{n^3}{\alpha} \right) \EE[\Vol_n(K)].
$$
The upper bound is established, and we have proven Theorem \ref{thmapprox}.

\section{Approximate scaling invariance of the facet distribution}

In this section we will prove Theorem \ref{thmscaling}.  \\ \\
We begin by fixing some compact and non-degenerate family of simplices $\mathcal{C} \subset \mathcal{S}$. Recall the definition of the measure $q(\cdot)$ (equation \ref{defq}). We have
\begin{equation} \label{eqMKC}
M_K(\mathcal{C}) = \EE \left [\int_{\Delta_n} \mathbf{1}_{ \{ F_r \in \mathcal{C} \} } dq(r) \right ].
\end{equation}
Our main ingredient will be the following lemma which helps us express the behavior of facets of $K$ in terms of limits of the expected behavior of facets of $K_\alpha$. This, in turn, will allow us to use the machinery developed in section 3 in order to calculate the right hand side of the above equation.

\begin{lemma} \label{toughlemma}
Let $f: \Delta_n \to [0,\infty)$ be a random function satisfying the following conditions: \\
(i) For all $r \in \Delta_n$, $f(r)$ is measurable with respect to the $\sigma$-algebra generated by the Brownian motion $B(\cdot)$. \\
(ii) There exist constants $C>0$ and $p \geq 1$ such that almost surely,
\begin{equation} \label{fcompV}
f(r) \leq C V(r)^p, ~~ \forall r \in \Delta_n
\end{equation}
where $V(r) = Vol_{n-1}(F_r)$ is defined in section 3. \\
(iii) The function $f$ is almost surely continuous in $\Delta_n$. \\
Then we have
$$
\EE \left [ \int_{\Delta_n} f(r) d q(r) \right ] = \lim_{\alpha \to \infty} \alpha^n \int_{\Delta_n} \EE \left [ f(r) \mathbf{1}_{E_\alpha(r)}  \right ] dr.
$$
\end{lemma}
We postpone the proof of this lemma to the end of the section. \\ \\

This lemma encourages us to define for all $\delta > 0$
$$
f_\delta(r) = \max \left (0, 1 - \delta^{-1} d_H(F_r, \mathcal{C})  \right ).
$$
where $d_H(s, \mathcal{C})$ denotes the minimal Hausdorff distance between $s$ and some $s' \in \mathcal{C}$. The non-degeneracy assumption of the family $\mathcal{C}$ implies
$$
F_r \in \mathcal{C} \Rightarrow V(r) \geq c
$$
for some constant $c>0$. Therefore, since the volume of a simplex is continuous with respect to the Hausdorff metric and by the compactness of $\mathcal{C}$, we deduce that there exists some $\delta_0,C>0$ such that 
$$
f_\delta(r) \leq C V(r), ~~ \forall 0 < \delta < \delta_0.
$$
Consequently, the assumption (\ref{fcompV}) is satisfied and we may apply the above lemma to get
$$
\EE \left [ \int_{\Delta_n} f_\delta(r) d q(r) \right ] = \lim_{\alpha \to \infty} \alpha^n \int_{\Delta_n} \EE  \left [ f_\delta(r) \mathbf{1}_{E_\alpha(r)}  \right ] dr.
$$
Since $f_\delta(r)$ is increasing with respect to $\delta$, the monotone convergence theorem and equation (\ref{eqMKC}) teach us that
$$
\lim_{\delta \to 0^+} \EE \left [ \int_{\Delta_n} f_\delta(r) d q(r) \right ] = M_K(\mathcal{C}).
$$
The last two equations give
$$
M_K(\mathcal{C}) = \lim_{\delta \to 0^+} \lim_{\alpha \to \infty} \alpha^n \int_{\Delta_n} \EE \left [ f_\delta(r) \mathbf{1}_{E_\alpha(r)}  \right ] dr.
$$
Recall that the equivalence class of $F_r$ and the event $E_\alpha(r)$ are independent (by Lemma \ref{lemfacets}). It follows that
\begin{equation} \label{eq444}
M_K(\mathcal{C}) = \lim_{\delta \to 0^+} \lim_{\alpha \to \infty} \alpha^n \int_{\Delta_n} p_\alpha(r) \EE \left [ \max \left (0, 1 - \delta^{-1} d_H(F_r, \mathcal{C})  \right ) \right ] dr.
\end{equation}
Let $\Gamma_1,...,\Gamma_{n-1}$ be independent standard Gaussian random vectors. For a point $r \in \Delta_n$, we define
\begin{equation} \label{defxr}
X_r = \Conv \left ( \{0\} \cup \left (\bigcup_{i=1}^{n-1} \left \{ \sum_{j=1}^i \sqrt{r_{j+1}} \Gamma_{j} \right \} \right ) \right ),
\end{equation}
so $X_r$ has the same distribution as $F_r$, up to a translation. Using formula (\ref{facets2}), equation (\ref{eq444}) becomes
$$
M_{K}(\mathcal{C}) = \lim_{\delta \to 0^+} \lim_{\alpha \to \infty} \alpha^n \int_{\Delta_n} p_\alpha(r) g(\delta, r) dr = 
$$
$$
\lim_{\delta \to 0^+} \lim_{\alpha \to \infty} \frac{2}{\pi} \int_{\Delta_n} \left (\prod_{j=2}^{n} \frac{1 - e^{-\alpha r_j}}{r_j} \right ) \frac{1}{\sqrt{r_1 r_{n+1}}}  \Phi(\alpha r_1) \Phi(\alpha r_{n+1}) g(\delta, r) dr
$$
where
$$
g(\delta, r) := \EE \left [ \max \left (0, 1 - \delta^{-1} d_H(X_r, \mathcal{C})  \right ) \right ].
$$
The bound (\ref{estphi}) suggests that the dominated convergence theorem may be used to attain
$$
M_{K}(\mathcal{C}) = \lim_{\delta \to 0^+} \frac{2}{\pi} \int_{\Delta_n} \left (\prod_{j=2}^{n} \frac{1}{r_j} \right )  \frac{1}{\sqrt{r_1 r_{n+1}}} g(\delta, r) dr.
$$
Moreover, we have
$$
\lim_{\delta \to 0^+} \max \left (0, 1 - \delta^{-1} d_H(X_r, \mathcal{C})  \right ) = \mathbf{1}_{ \{ X_r \in \mathcal{C} \}}.
$$
an application of the monotone convergence theorem with the two previous equations gives
$$
M_{K}(\mathcal{C}) = \frac{2}{\pi} \int_{\Delta_n} \left (\prod_{j=2}^{n} \frac{1}{r_j} \right )  \frac{1}{\sqrt{r_1 r_{n+1}}} \PP(X_r \in \mathcal{C}) dr =
$$
$$
\frac{2}{\pi} \int_{ \Delta_{n-1} } \left ( \prod_{j=2}^n \frac{1}{r_j} \right ) \PP(X_r \in \mathcal{C}) \left ( \int_0^{1 - \sum_{i=2}^{n} r_i} \frac{1}{\sqrt{t(1 - \sum_{i=2}^{n} r_i - t)}} dt \right ) dr_2 ... dr_{n} = 
$$
$$
2 \int_{ \Delta_{n-1}  } \left ( \prod_{j=2}^n \frac{1}{r_j} \right ) \PP(X_r \in \mathcal{C}) dr_2 ... dr_n.
$$
Recall that for a set $\mathcal{C} \subset \mathcal{S}$ and $t > 0$, we understand $t \mathcal{C}$ as $\{t s; ~ s \in \mathcal{C} \}$. 
Note that $X_{tr} \sim \sqrt{t} X_{r}$. Therefore, by substituting $r_i = \lambda w_i$ for $2 \leq i \leq n$ in the last integral we attain
that, for all $\lambda \geq 1$,
$$
M_K(\mathcal{C}) = 2 \int_{ \Delta_{n-1} } \left ( \prod_{j=2}^n \frac{1}{r_j} \right ) \PP(X_r \in \mathcal{C}) dr_2 ... dr_n \leq
$$
$$
2 \int_{ \lambda \Delta_{n-1} } \left ( \prod_{j=2}^n \frac{1}{r_j} \right ) \PP(X_r \in \mathcal{C}) dr_2 ... dr_n =
$$
$$
2 \int_{ \Delta_{n-1} } \left ( \prod_{j=2}^n \frac{1}{w_j} \right ) \PP(X_{\lambda w} \in \mathcal{C}) dw_2...dw_n = M_K \left (  \frac{1}{\sqrt{\lambda}} \mathcal{C} \right ).
$$
It follows that $M_K(t \mathcal{C})$  is a decreasing function of $t$, which means that
$\lim_{t \to 0} M_K(t \mathcal{C})$ exists in the wide sense. This completes the first part of the theorem. \\ \\

The second part of the theorem will follow directly from the next technical lemma.
\begin{lemma}
Let $\mathcal{C} \subset \mathcal{S}$ be compact in the Hausdorff metric, such that every $s \in \mathcal{C}$
has a non-empty relative interior. Then
$$
\int_{ \RR_+^{n-1} } \left ( \prod_{j=2}^n \frac{1}{r_j} \right ) \PP(X_r \in \mathcal{C}) dr_2... dr_n < \infty,
$$
where $\RR_+^{n-1} = \{(r_2,...,r_n); ~ r_i \geq 0, ~~ \forall 2 \leq i \leq n \}$.
\end{lemma}
\emph{Proof:}
For a simplex $s \in \mathcal{C}$, let $M(s)$ and $m(s)$ denote the length of the longest and shortest two-dimensional
edge of $s$, respectively. Since $\mathcal{C}$ is compact we may define
$$
M = \max_{s \in \mathcal{C}} M(s), ~~ m = \min_{s \in \mathcal{C}} m(s).
$$
Along with the assumption that the simplices are non-degenerate (e.g. have a non-empty interior), the reader can easily verify that $m > 0$.
Fix  $(r_1,...,r_n) = r \in \RR_+^n$ and for all $1 \leq i \leq n$ let $t_i = \sum_{j=1}^i r_i$. 
Let $\Gamma$ be a standard Gaussian random vector in $\RR^n$. Recall the definition of $X_r$ in equation (\ref{defxr}). We estimate
$$
\PP(m(X_{r}) > m/2) \leq \PP \left (\min_{2 \leq i \leq n}  |B(t_{i-1}) - B(t_{i})| > m/2  \right) \leq
$$
$$
\PP \left ( \sqrt{\min_{2 \leq i \leq n} r_i} |\Gamma| > m/2 \right ) \leq C_1 \exp \left (-c_1 \frac{m^2}{\min_{2 \leq i \leq n} r_i} \right ) \leq
$$
$$
C_1 \exp \left (-c_1 m^2 \prod_{i=2}^n r_i^{-1/{n-1}} \right )
$$
for some constants $c_1,C_1>0$. On the other hand, we can write
$$
\PP(M(X_{r}) < 2M) \leq \PP( \max_{2 \leq i \leq n} \sqrt{r_i} |\Gamma| < 2M) \leq 
$$
$$
C_2 \max_{2 \leq i \leq n} r_i^{-n/2} \leq C_2 \prod_{i=2}^n r_i^{-1/2}.
$$
for some constants $c_2, C_2>0$ (which may depend on $M$ and $n$). Using these two estimates, we obtain
$$
\PP(M(X_{r}) \in \mathcal{C}) \leq C \min \left (\prod_{i=2^n} r_i^{-1/2}, \exp \left (- c \prod_{i=2}^n r_i^{-1/{n-1}} \right ) \right ) 
$$
for some constants $c,C>0$. We therefore get
$$
\int_{ \RR_+^{n-1} } \left ( \prod_{j=2}^n \frac{1}{r_j} \right ) \PP(X_r \in \mathcal{C}) dr_2...dr_n \leq 
$$
$$
\int_{ \RR_+^{n-1} } C \min \left (  \prod_{j=2}^n \frac{1}{r_j^{3/2}} , \prod_{j=2}^n \frac{1}{r_j} \exp \left (- c \prod_{i=2}^n r_i^{-1/{n-1}} \right )  \right ) dr_2...dr_n \leq .
$$
$$
\int_{ \RR_+^{n-1} } \min \left (  \prod_{j=2}^n \frac{1}{r_j^{3/2}} , C' \right ) dr_2...dr_n
$$
for some constant $C'>0$. The above integral obviously converges, and the lemma is proven. \qed \\

\begin{remark}
The method above may actually be used to find a precise formula for, say, the expected number of facets of $K$ whose volume is between two given constants, in terms of the distribution function of the product of $\chi$-random variables. As another example, by taking $f(r) = V(r)^2$ in Lemma \ref{toughlemma}, another quantity which we can easily calculate using this method is the expected volume of the facet containing a point $x \in \partial(K)$ when $x$ is uniformly generated on the set $\partial(K)$. Quantities of this sort may serve us to attain a little more information on the distribution of facets of $K$.
\end{remark}
$$
~
$$
\textbf{Proof of Lemma \ref{toughlemma}:} \\
We divide the proof into several steps. \\
\textbf{Step 1.} 
We start with showing that for all $A \subset Int(\Delta_n)$ compact, one has almost surely
\begin{equation} \label{eqstep1}
q(A) = 0 \Rightarrow \limsup_{\alpha \to \infty} q_\alpha(A) = 0. 
\end{equation}
Clearly, it suffices to fix an arbitrary continuous function $B(t)$ and show that the last equation holds almost surely with respect to the Poisson process. Set $A \subset Int(\Delta_n)$ and assume that $q(A) = 0$ (note that the last event is measurable with respect to the Brownian motion $B(\cdot)$. Its only property that we will use is continuity). For any $r \in A$, denote
$$
T(r) = \bigl \{t \in [0,1] \mid ~ \langle n_r,  B(t) \rangle > H(r)  \bigr \}.
$$
The assumption $q(A) = 0$ implies that $T(r) \neq \emptyset$ for all $r \in A$.  By the continuity of the Brownian motion, it follows that $T(r)$ is an open set and moreover $\lambda(T(r))$ (namely, the Lebesgue measure of $T(r)$) is a positive and continuous function and so attains a minimum on $A$, which we denote by $\epsilon$. Now, for every $r \in A$, let $\tilde T(r)$ be some closed subset of $T(r)$ whose Lebesgue measure is at least $\epsilon / 2$. It is easy to check that for all $r \in \Delta_n$ there exists an open neighborhood, $N(r)$ such that 
$$
r' \in N(r) \Rightarrow \tilde T(r) \subset T(r').
$$
Since $A$ is compact, there exists a finite set $r_1,...,r_N \in A$ such that 
$$
A \subset N(r_1) \cup ... \cup N(r_N).$$
Now, suppose  that for each $1 \leq i \leq N$, one has $\Lambda_\alpha \cap \tilde T(r_i) \neq \emptyset$. Then for all $r \in A$ we have $T(r) \supset \tilde T(r_j)$ for some $1 \leq j \leq N$, which implies that $T(r) \cap \Lambda_\alpha \neq \emptyset$, which in turn means that $E_\alpha(r)$ does not hold. It follows that
$$
\PP(q_\alpha(A) = 0 | B(t)) \geq \PP \left ( \bigcap_{1 \leq j \leq N} \{ \Lambda_\alpha \cap \tilde T_i \neq \emptyset  \}  \right ) \geq 1 - N e^{-\epsilon \alpha / 2}.
$$
Since $\Lambda_\alpha$ is monotone in $\alpha$ (with respect to inclusion), it follows that the events $\{ \Lambda_\alpha \cap \tilde T_i \neq \emptyset \}$ are increasing with $\alpha$, and we get $\limsup_{\alpha \to \infty} q_\alpha(A) = 0$ almost surely. \\ \\
\textbf{Step 2} We show that if $A \subset Int(\Delta_n)$ is compact and $f$ is supported in $A$ then, almost surely,
\begin{equation} \label{eqstep2}
\limsup_{\alpha \to \infty} \int_{\Delta_n} f(r) d q_\alpha(r) \leq \int_{\Delta_n} f(r) dq(r).
\end{equation}
Again, we may fix an arbitrary continuous function $B(\cdot)$ and prove that the last equation holds almost surely with respect to the Poisson process. Observe that for almost every point $r \in \Delta$, there exists a neighborhood $M(r)$ (which depends on $B(t)$) such that 
$$
r_1, r_2 \in M(r) \Rightarrow \mathrm{Int}(\Conv(F_{r_1}, 0)) \cap \mathrm{Int}(\Conv(F_{r_2}, 0)) \neq \emptyset.
$$
We will use the following geometric fact: for any polytope $P$ which contains the origin, for any two distinct facets $F_1,F_2 \subset \partial P$ (of co-dimension $1$) one has $\mathrm{Int}(\Conv(F_{1}, 0)) \cap \mathrm{Int}(\Conv(F_{2}, 0)) = \emptyset$. Since $0 \in K_\alpha$ by definition, this fact implies that almost surely
\begin{equation} \label{onepoint}
\limsup_{\alpha \to \infty} q_\alpha(M(r)) \leq 1.
\end{equation}
Denote $S = supp(q) \cap A$. For all $\epsilon > 0$, define 
$$
M_\epsilon = \bigcup_{r \in S} \left ( M(r) \cap B(r, \epsilon) \right )
$$
where $B(r, \epsilon)$ denotes the open Euclidean ball of radius $\epsilon$ centered at $r$. Then due to (\ref{onepoint}), we learn that almost surely
$$
\limsup_{\alpha \to \infty} \int_{M_\epsilon} f(r) d q_\alpha(r) \leq \sum_{r \in S} \sup_{r' \in B(r,\epsilon)} f(r').
$$
Next, since $M$ is open by definition, using (\ref{eqstep1}), we have
$$
\limsup_{\alpha \to \infty} \int_{A \setminus M_\epsilon} f(r) d q_\alpha(r) = 0.
$$
Evidently, the set function $\limsup_{\alpha \to \infty} q_\alpha(\cdot)$ is sub-additive. Thus, the last two equations imply that
for all $\epsilon > 0$ we have almost surely
$$
\limsup_{\alpha \to \infty} \int_{\Delta_n} f(r) d q_\alpha(r) \leq \sum_{r \in S} \sup_{r' \in B(r,\epsilon)} f(r').
$$
Recall that $f$ is assumed to be continuous and compactly supported. By using the uniform continuity of $f$ and taking $\epsilon \to 0$ we attain (\ref{eqstep2}). \\ \\
\textbf{Step 3:} 
We show that almost surely,
\begin{equation} \label{eqstep3}
\lim_{\alpha \to \infty} \int_{\Delta_n} f(r) d q_\alpha(r) = \int_{\Delta_n} f(r) d q(r).
\end{equation}
Assume by contradiction that one has with positive probability $\epsilon > 0$, where
\begin{equation} \label{contradiction}
\epsilon = \int_{\Delta_n} f(r) d q(r) - \liminf_{\alpha \to \infty} \int_{\Delta_n} f(r) d q_\alpha(r).
\end{equation}
Recall that $V(r) = \Vol_{n-1} (F_r)$. By Corollary \ref{corlimit}, we have almost surely
\begin{equation} \label{intislarge}
\liminf_{\alpha \to \infty} \int_{\Delta_n} V(r) d q_\alpha(r) \geq \int_{\Delta_n} V(r) d q(r).
\end{equation}
Now, according to the assumption (\ref{fcompV}) we have for all $r \in \Delta_n$,
$$
f(r) \leq C V(r)^p \leq C V(r) \max_{r \in \Delta_n} V(r)^{p-1}.
$$
Consequently, almost surely there exists a constant $m>0$ such that the function $h(r) = V(r) - \frac{m}{2} f(r)$ satisfies $h(r) \geq 0$ for all $r \in \Delta_n$. Equations (\ref{intislarge}) and (\ref{contradiction}) imply,
\begin{equation}
\limsup_{\alpha \to \infty} \int_{\Delta_n} h(r) dq_\alpha(r) \geq \int_{\Delta_n} h(r) dq + \frac{m \epsilon}{2}.
\end{equation}
Since $h(r)$ is continuous, there exists a compact $A \subset Int(\Delta_n)$ such that
\begin{equation}
\limsup_{\alpha \to \infty} \int_{A} h(r) dq_\alpha(r) \geq \int_{A} h(r) dq + \frac{m \epsilon}{4}
\end{equation}
with positive probability. This contradicts (\ref{eqstep2}). Therefore, we have established (\ref{eqstep3}). \\ \\
\textbf{Step 4} 
To finish the proof of the lemma, we argue that the dominated convergence theorem may be used to show that
\begin{equation} \label{domiconv}
\EE \left[ \lim_{\alpha \to \infty} \int_{\Delta_n} f(r) dq_\alpha(r) \right ] = \lim_{\alpha \to \infty} \EE \left[ \int_{\Delta_n} f(r) dq_\alpha(r) \right ].
\end{equation}
Indeed, according to (\ref{volsa1}) and by the assumption (\ref{fcompV}), we have
$$
\int_{\Delta_n} f(r) dq_\alpha(r) \leq C \int_{\Delta_n} V(r)^p dq_\alpha(r) \leq
$$
$$
C'_{n,p} C \left (\int_{\Delta_n} V(r) dq_\alpha(r) \right )^p = C'_{n,p} C \Vol_{n-1} (\partial K_\alpha)^p
$$
where in the second inequality we use the fact that $q_\alpha$ is a counting measure, and $C'_{n,p}$ is a positive constant depending only on $n$ and on $p$. The inclusion 
$$
K_\alpha \subset K \subset \left  \{x \in \RR^n; ~ |x| \leq \max_{0 \leq t \leq 1} |B(t)| \right \}
$$ 
teaches us that
$$
\Vol_{n-1} (\partial K_\alpha)^p \leq \bigl ( \Vol_{n-1} (\partial \BBN) \bigr)^p \left (\max_{0 \leq t \leq 1} |B(t)| \right )^{p(n-1)},
$$
and therefore
$$
\int_{\Delta_n} f(r) dq_\alpha(r) \leq C \left (\max_{0 \leq t \leq 1} |B(t)| \right )^{p(n-1)}
$$
where $C$ is deterministic. Now, it is well known that the maximum of a Brownian motion in a bounded interval has finite moments of all orders (for example, by a combination of Doob's theorem with the fact that the Gaussian distribution has finite moments). In other words,
$$
\EE \left [ \left (\max_{0 \leq t \leq 1} |B(t)| \right )^{p(n-1)} \right ] < \infty.
$$
The two last inequalities and the dominated convergence theorem finally prove (\ref{domiconv}). We now combine (\ref{domiconv}), (\ref{eqstep3}) and (\ref{superfubini2}) to get
$$
\EE \left [ \int_{\Delta_n} f(r) d q(r) \right ] = \EE \left[ \lim_{\alpha \to \infty} \int_{\Delta_n} f(r) dq_\alpha(r) \right ] = 
$$
$$
\lim_{\alpha \to \infty} \EE \left[ \int_{\Delta_n} f(r) dq_\alpha(r) \right ] =  \lim_{\alpha \to \infty} \alpha^n \int_{\Delta_n} \EE \left [ f(r) \mathbf{1}_{E_\alpha(r)}  \right] dr.
$$
The lemma is complete.
\qed \\ \\

We finish this section with a small remark. As a consequence of the above lemma, and by Corollary \ref{corlimit}, we have
$$
\EE \left [ \int_{\Delta_n} V(r) d q(r) \right ] =
$$
$$
\lim_{\alpha \to \infty} \alpha^n \int_{\Delta_n} \EE \left [ V(r) \mathbf{1}_{E_\alpha(r)} \right ] dr = 
$$
$$
\lim_{\alpha \to \infty} \EE \left [ \Vol_{n-1}(\partial K_\alpha) \right ]  =  \EE \left [ \Vol_{n-1}(\partial K) \right ].
$$
Note that a-priori, the quantity $\int_{\Delta_n} V(r) d q(r)$ need not be equal to the surface area of $K$, because some of its surface area may be contained in facets of dimension smaller than $n-1$. We therefore have the following corollary to the above lemma:

\begin{corollary} \label{corfacets}
Let $A \subset \partial K$ be the set of points not contained in the interior of any $(n-1)$-dimensional facet of $K$. We have
$$
\Vol_{n-1} (A) = 0.
$$
\end{corollary}

\section{Comments and Possible Further Research} \label{secdiscussion}

\subsection{Higher moments}
In this paper, we derived formulas for the expectation (i.e, the first moment) of the volume and surface area of $K$. It may be interesting to 
investigate the behaviour of higher moments. In particular, it may be interesting to ask, for example, how \emph{concentrated} the volume of $K$ is around its mean. \\ \\

One possible strategy of finding the second moment is by using an analogous formula to (\ref{volcalc2}). It can be seen that, in the notations of section 3,
$$
\EE[\Vol_n(K)^2] = 
$$
$$
\frac{1}{n^2} \int_{\Delta_n \times \Delta_n}  \EE \left [ V(x) V(y) R(x) R(y) \mathbf{1}_{E_\alpha(y)} \mathbf{1}_{E_\alpha(x)} \right ]  dx dy.
$$
It seems rather hard to find a precise formula for the integrand. However, when the dimension $n$ is large, it is possible to show using standard techniques related to high-dimensional measure concentration that for a typical choice of $x,y \in \Delta_n \times \Delta_n$ (i.e, with high probability), the facets $F_x$ and $F_y$ will be approximately orthogonal to each other. Now, note that the random variables
$R(x)$ and $\mathbf{1}_{E_\alpha(x)}$ depend only on $\langle B(t), n_x \rangle$ while the variables $R(y)$ and $\mathbf{1}_{E_\alpha(y)}$ depend only on $\langle B(t), n_y \rangle$. It follows that when $n_x \perp n_y$, these two pairs are mutually independent. Since $n_x$ and $n_y$ are \emph{almost} orthogonal, it is reasonable to suspect that for a typical pair $x,y$ one would have
$$
\EE \left [ V(x) V(y) R(x) R(y) \mathbf{1}_{E_\alpha(y)} \mathbf{1}_{E_\alpha(x)} \right ] \approx
$$
$$
\EE \left [ V(y) R(y) \mathbf{1}_{E_\alpha(y)}  \right ] 
\EE \left [ V(x) R(x) \mathbf{1}_{E_\alpha(x)}  \right ].
$$
In some sense, the above would imply that $\EE[\Vol_n(K)]^2$ is close to $\EE[\Vol_n(K)^2]$ when the dimension is large. This suggests that the answer to the following question may be positive:

\begin{question}
Is it true that
$$
\lim_{n \to \infty} \frac{\sqrt{Var[\Vol_n(K)]} }{\EE[\Vol_n(K)]} = 0?
$$
\end{question}

\subsection{Smoothness}

In 1983, El Bachir (\cite{Elb}) proved the assertion of P. L\'{e}vy that almost surely, the convex hull of a planar Brownian motion has a smooth boundary. Later, in 1989, Cranston, Hsu and March (\cite{CHM}) showed that in fact, it is exactly H\"{o}lder $1 \frac{1}{2}$-smooth. A natural question would be about the extension of these facts to higher dimensions:

\begin{question}
Does the body $K$ have a smooth boundary for any dimension? 
\end{question}

 It doesn't seem straightforward to adapt their methods even to the three dimensional case: their proofs rely on the fact that if a $2$-dimensional convex hull has a "corner", then the directions of its supporting hyperplanes will contain some interval, which will in turn contain a rational direction. This fact allows them to express the smoothness of the boundary as a countable intersection of events, each depending on a single direction. It is easy to see that in $3$-dimensional space, this is already not the case: it is not hard to construct a body whose boundary is not smooth, but a uniformly generated random $2$-dimensional projection of this body will almost surely be smooth. \\ \\

The following heuristic argument may suggest that the boundary is, in fact, smooth in higher dimensions:  \\ \\
For two $(n-1)$-dimensional facets $s,t$ of a polytope $P$, we say that $s,t$ are neighbors if the intersection of $s,t$ has Hausdorff dimension $n-2$. The first step is to try to prove that for any $\epsilon > 0$, the probability that there exists at least one pair of $n-1$-dimensional neighboring facets of $K_\alpha$ such that the angle between the two is at least $\epsilon$ goes to zero as $\alpha \to \infty$. The idea is the following: in the notation of section 3, note that any choice of two neighboring facets of $K_\alpha$ corresponds to a choise of $n+1$ points from the process $\Lambda_\alpha$. In other words, it corresponds to the choice of two point $r, s \in \Delta_n$ such that $t(r)$ and $t(s)$ differ by only one coordinate. \\ \\
Consider the event $E$ that both $F_{r}$ and $F_{s}$ are facets of $K_\alpha$ and that the angle between the two is more than $\epsilon$. By the representation theorem for the Brownian bridge, it follows that one can first generate the points $B(t)$ for $t \in t(r) \cup t(s)$, and then "fill in" the missing gaps by Brownian bridges, as carried out in the proof of Lemma \ref{lemfacets}. Now, project the Brownian motion onto the two dimensional subspace spanned by $n_{r}$ and $n_{s}$. Following the same lines as the proof of Lemma \ref{lemfacets}, we see that the event $E$ is reduced to the event that $n-1$ independent discrete Brownian bridges and two Brownian motions all stay in a wedge of angle $\pi - \epsilon$ (in Lemma \ref{lemfacets}, the event $E_{\alpha}(r)$ was equivalent to the fact that they stay in a half-space). \\ \\
The next step would be to generalize the bounds in Lemma \ref{lemsec2} to a wedge rather than a half-space. Considering the bounds obtained in \cite{D} as well as by the conformal invariance of the Brownian motion, it is reasonable to expect that the probability of an $\alpha$-step random walk to stay in a wedge of angle $\theta$ is of the order $\alpha^{- \frac{\pi}{2 \theta}}$, and that for an $\alpha$-step discrete Brownian bridge to stay in such a wedge, would be of order $\alpha^{- \frac{\pi}{\theta}}$. If these estimates are indeed correct, plugging them in the analogue of equation (\ref{facets2}) for wedges should give roughly
$$
\PP(E) \leq \alpha^{-1} \prod_{i=2}^n (\alpha r_i)^{- \frac{\pi}{\theta}} (\alpha^2 r_1 r_{n+1})^{- \frac{\pi}{2 \theta}}.
$$
If the above bound is true, it would imply that
$$
\alpha^{n+1} \int_{\Delta_{n+1}} \PP(E) d(r,s) \to 0
$$
as $\alpha \to \infty$. This means that when $\alpha$ is large, one should expect $K_\alpha$ to be "smooth" in the sense that any neighboring two faces have angle less than $\epsilon$ between them. Next, an analogue of Lemma \ref{toughlemma} may be used to show that this property is preserved when passing to the limit: Indeed, the lemma shows that any facet of $K$ already becomes a facet of $K_\alpha$ for $\alpha$ large enough. Therefore, if $K$ has two neighboring facets with some angle $\epsilon > 0$ between them, then $K_\alpha$ will have two such facets for all $\alpha$ large enough, and we would arrive at a contradiction.

\subsection{Neighborliness of the approximating polytope}

A polytope $P \subset \RR^n$ is said to be \emph{$k$-neighborly} if for any choice of $k$ vertices of $P$, $v_1,...,v_k$, the simplex
$\Conv(v_1,...,v_k)$ is a facet of $P$. The concept of neighborlyness is related to the ability of linear programming to find solutions to systems of underdetermined linear equations (see \cite{DT}). \\ \\

We believe that a slight generalization of the method introduced in \cite[section 2]{E1} may be used 
to show that when $\alpha$ is a polynomial of $n$, the polytope $K_\alpha$ is $(c n / \log^2 n)$-neighborly, with probability approaching $1$ as $n \to \infty$. The idea of proof is as follows: \\ \\
Fix some $k<n$ and take $\alpha = n^{10}$, say. Let $0 \leq s_1 <... < s_k \leq 1$ be a selection of $k$ points from $\Lambda_\alpha$. Define $F = \Conv(B(s_1),...,B(s_k))$. Let us try to understand the event that $F$ is contained in the boundary of $K_\alpha$. By the representation theorem of the Brownian bridge, one may first generate the points $B(s_1),...,B(s_k)$ and then "fill in" the gaps with Brownian bridges, as carried out in the proof of Lemma \ref{lemfacets}. In view of this, and by considering the projection of the Brownian motion onto $F^\perp$, the event above is reduced to the fact that $k-1$ discrete Brownian bridges of length smaller than $n^{10}$ in $\RR^{n-k}$ and two random walks are all contained in some open halfspace of $\RR^{n-k}$. At this point, Theorem 2.1 in \cite{E1} comes to our service. According to this theorem, for a random walk of polynomial length in $\RR^n$ with probability $1 - n^{-10}$, there exists a unit vector whose scalar product with any internal point of the random walk is of order $\log n$. By proving an analogous result for a discrete Brownian bridge, one would be able to combine $\log n$ such vectors together to create a vector separating these bridges from the origin, thus proving that $F$ is in the boundary of $K_\alpha$. \\ \\

In some sense, the property of a polytope being $k$-neighborly is contradictory to the fact that it approximates a smooth convex body: it is not hard to realize that as a family of polytopes approaches a smooth convex body, they become less neighborly in some sense. Therefore, the above fact may be interesting considering the fact that $K_{n^4}$ is already a good volumetric approximation for the polytope $K$ which is obviously not $2$-neighborly (and may, in fact, be smooth). \\ \\

\subsection{Comparison with a Gaussian Polytope}
Fix a dimension $n$. For an integer $\alpha \in \mathbb{N}$, let $\Gamma_1,...,\Gamma_\alpha$ be independent standard Gaussian random vectors in $\RR^n$.
We define $G_\alpha = \Conv(\Gamma_1,...,\Gamma_\alpha)$, the convex hull of independent Gaussian points. The object $G_\alpha$ is usually referred to as a \emph{Gaussian polytope}. The study of the Gaussian polytope began in the 60's by R\'enyi and Sulanke, and ever since it has been deeply investigated (see \cite{as,hr} and \cite{CMR}). \\

For geometers in convexity, a central motivation in the investigation of random polytopes stems from the fact that random objects often admit a rather pathological behavior and thus often serve as counterexamples to certain conjectures. A-priori, one may have expected that the Gaussian polytope may be used to serve as a counter example for certain phenomena related to the distribution of mass on convex bodies. Alas, results such as \cite{kk} and \cite{fleury} suggest that this object admits a quite regular and symmetric nature. It may therefore be interesting to try to find a construction of random polytopes that are, in some sense, as irregular and asymmetric as possible. Some of the estimate we obtained here may point to the fact that the behavior of the polytope $K_\alpha$ is less regular than the behavior of $G_\alpha$. We list some possible qualitative differences between those two constructions. \\ \\

It can be seen (see e.g., \cite{fleury}, Theorem 3) that the facets of $G_n$ admit a rather regular behavior in the sense that for any
given family of simplices $\mathcal{C}$, the function $M_{G_\alpha}(t \mathcal{C})$ is rather concentrated around a specific value of $t$. In other words, there is a typical "correct" scale at which one expects to find most of the facets of $G_\alpha$. On the other hand, by Theorem \ref{thmscaling}, we know that when $\alpha$ is rather large, one expects to find facets of the same shape at a wide range of scales, which suggests a less regular behavior. \\ \\

As shown in \cite{kk}, the covariance matrix of a uniform point randomly generated of $G_\alpha$ is not far from a multiple of the identity matrix (for $\alpha$ large enough). This property is sometimes referred to as \emph{isotropicity}. It can be also seen by \cite[Theorem 3]{fleury} that the covariance matrix of a typical facet of $G_\alpha$ is rather isotropic. On the other hand, in view of the formulas developed in sections 3 and 6, one may expect that the polytope $K_\alpha$ exhibits a very different behavior. Formula (\ref{facets}), suggests that covariance matrix of a uniform point on $K_\alpha$, as well as on one of the large-scale facets of $K_\alpha$ will have a covariance matrix close to the one of $B(\cdot)$, which is far from isotropic. Note also that the small-scale facets which, in the notations of section 3, have $r_i \ll 1$ for $2 \leq i \leq n$, seem to "ignore" the covariance structure of $B(\cdot)$. Therefore, it is reasonable to guess that there is a phase shift in the geometry of facets of $K_\alpha$ when passing from larger to smaller scales of facets. \\ \\

Both constructions, $G_\alpha$ and $K_\alpha$ seem to tend to a rather smooth shape as $\alpha \to \infty$: it is well known that the shape of $G_\alpha$ becomes quite close to a Euclidean ball when the value of $\alpha$ is large enough, namely exponential in $n$, and $K_\alpha$ approaches $K$, which as the last section suggests, may have a smooth boundary. It is known that $G_\alpha$ is a highly-neighborly polytope when $n$ is a proportion of $\alpha$ (see \cite{DT}), but it becomes almost non-neighborly as it approaches a Euclidean ball. On the other hand, $K_{n^{10}}$ which is, by Theorem \ref{thmapprox}, already close in expectation to its "smooth limit", may be a highly neighborly polytope, as the previous subsection suggests. \\ \\ \\

\emph{Acknowledgements}
I would like to thank Itai Benjamini for fruitful discussions. I am also very thankful to Christoph Th\"{a}le for pointing him to the observation that Theorem 1 may be generalized to obtain formulae for all intrinsic volumes. Finally, I thank the two anonymous referees of this paper for doing a very thorough job reviewing it, providing me with numerous helpful and enlightening comments.

\bigskip {\noindent Department of Mathematics, Weizmann Institute of Science, Rehovot, Israel \\  {\it e-mail address:}
\verb"roneneldan@gmail.com" } 

\end{document}